\documentclass{article}
\usepackage{latexsym,amsmath,amssymb,amscd}
\usepackage[dvips]{color}
\begin{document}

\title{Note on the descriptions of the Euler-Poisson equations in various co-ordinate systems }
\author{Tetu Makino \footnote{Professor Emeritus at Yamaguchi University, Japan. e-mail: makino@yamaguchi-u.ac.jp}}
\date{\today}
\maketitle

\newtheorem{Lemma}{Lemma}
\newtheorem{Proposition}{Proposition}
\newtheorem{Theorem}{Theorem}
\newtheorem{Definition}{Definition}
\newtheorem{Remark}{Remark}
\newtheorem{Formula}{Formula}
\numberwithin{equation}{section}

\begin{abstract}
In this note we derive the descriptions of the system of Euler-Poisson equations which governs the hydrodynamic evolution of gaseous stars
in various co-ordinate systems. This note does not contain essentially new 
results for astrophysicists, but mathematically rigorous derivations cannot be found in the literatures written by physicists so that it will be useful to prepare details of rather stupidly honest derivations of the equations in various flames when we are going to push
forward the mathematical research of the problem. \\



\end{abstract}

\section{Euler-Poisson equations}

\subsection{Euler-Poisson equations and Newton potential}

The Euler-Poisson equations which govern evolutions of a gaseous star are
\begin{subequations}
\begin{align}
&\frac{\partial\rho}{\partial t }+(\nabla |\rho\vec{v})=0 \label{1a}\\
&\rho\Big[\frac{\partial\vec{v}}{\partial t}+
(\vec{v}|\nabla)\vec{v}\Big]+\nabla P=-\rho\nabla \Phi \label{1b}\\
&\triangle\Phi =4\pi \mathsf{G}\rho. \label{1c}
\end{align}
\end{subequations}
The independent variable is $(t,x)=(t, x^1,x^2, x^3) \in [0,T)\times \mathbb{R}^3$.
$\mathsf{G}$ is a positive constant.
The unknown functions are  the density field $\rho=\rho(t,x)$, the pressure field $P=P(t,x)$, the gravitational potential $\Phi=\Phi(t,x)$, and
the velocity field $\vec{v}=(v^1,v^2,v^3)^{\top}(t,x)$. 

We are using the usual notations
\begin{align*}
&(\nabla| \rho\vec{v})=\sum_{k=1}^3\frac{\partial}{\partial x^k}
(\rho v^k), \\
&(\vec{v} | \nabla)v^j=\sum_{k=1}^3v^k\frac{\partial v^j}{\partial x^k}, 
\quad j=1,2,3,\\
&\nabla Q=\Big(\frac{\partial Q}{\partial x^1},
\frac{\partial Q}{\partial x^2},
\frac{\partial Q}{\partial x^3}\Big)\quad \mbox{for}\quad Q=P, \Phi \\
&\triangle \Phi=\sum_{k=1}^3\frac{\partial^2\Phi}{(\partial x^k)^2}.
\end{align*}\\

Also, we use the notation
\begin{equation}
\frac{D}{Dt}=\frac{\partial}{\partial t}+(\vec{v}|\nabla) 
=\frac{\partial}{\partial t}+\sum_{k=1}^3
v^k\frac{\partial}{\partial x^k},
\end{equation}\\
which rewrite \eqref{1a}, \eqref{1b} as
\begin{align*}
&\frac{D\rho}{Dt}+\rho (\nabla| \vec{v})=0, \\
&\rho \frac{D\vec{v}}{Dt}+\nabla P =-\rho\nabla\Phi.
\end{align*}\\

We put the following assumption:\\

{\bf (A): } {\it The pressure  $P$ is a given smooth function of $\rho>0$
such that $P>0, dP/d\rho >0$ for $\rho>0$ and
there is a smooth function $\Lambda \in C^{\infty}(\mathbb{R})$ such that $\Lambda(0)=0$ and
\begin{equation}
P=\mathsf{A}\rho^{\gamma}(1+\Lambda(\mathsf{A}\rho^{\gamma-1})),
\end{equation}
where $\mathsf{A}, \gamma$ are positive constants and 
$1<\gamma< 2$. }\\

\begin{Definition}
A solution $\rho=\rho(t,x), \vec{v}=\vec{v}(t,x),
\Phi=\Phi(t,x)$ will be called a {\bf compactly supported
classical solution} if $\rho, \vec{v},\Phi \in C^1([0,T)\times \mathbb{R}^3),
\Phi(t,\cdot)\in C^2(\mathbb{R}^3)$,
$\rho \geq 0$ everywhere, and the support of $\rho(t,\cdot)$
is compact for $\forall t\in [0,T)$. 
\end{Definition}
Without loss of generality we assume
$\vec{v}(t,\cdot)$ is bounded on $\mathbb{R}^3$, since any  modification of $\vec{v}$ outside
the support of $\rho$ is free.\\

%\subsection{}

For any compactly supported classical solution the Laplace equation 
(\ref{1c}) can be solved by the Newton potential
\begin{equation}
\Phi(t,x)=-\mathsf{G}\int
\frac{\rho(t,{x}')}{|\vec{x}-\vec{x}'|}
d\mathcal{V}({x}'). \label{3}
\end{equation}
Here $d\mathcal{V}({x})$ denotes the usual volume element
$dx^1dx^2dx^3$, and $$|\vec{x}|=\sqrt{(x^1)^2+(x^2)^2+(x^3)^2}$$ for $\vec{x}=(x^1,x^2,x^3)^T$.

In this note, we shall specify the solution $\Phi$ of (\ref{1c}) by this
Newton potential (\ref{3}) for any compactly supported classical solution.

\subsection{Conservation of mass, energy and angular momentum}

It is well known that the total mass
\begin{equation}
M:=\int\rho(t,x)d\mathcal{V}({x})
\end{equation}
and the total energy
\begin{align}
E&:=\int\Big(\frac{1}{2}\rho|\vec{v}|^2+\Psi(\rho)
+\frac{1}{2}\rho\Phi\Big)d\mathcal{V}({x}) \nonumber \\
&=\int\Big(\frac{1}{2}\rho|\vec{v}|^2+\Psi(\rho)\Big)d\mathcal{V}({x})
-\frac{\mathsf{G}}{2}
\int\int
\frac{\rho(t,x)\rho(t,{x}')}{|\vec{x}-\vec{x}'|}
d\mathcal{V}({x})d\mathcal{V}({x}')
\end{align}
are constants with respect to $t$ along
any compactly supported classical solution.

 Here the state quantities $u=u(\rho)$ (enthalpy) and
$\Psi(\rho)$ is defined by
\begin{equation}
u=\int_0^{\rho}\frac{dP}{\rho}, \qquad
\Psi(\rho)=\int_0^{\rho}ud\rho.
\end{equation}

\begin{Remark}
Note that $$
u=\frac{\mathsf{A}\gamma}{\gamma-1}\rho^{\gamma-1}, \qquad \Psi(\rho)=\frac{\mathsf{A}\rho^{\gamma}}{\gamma-1}=\frac{P}{\gamma-1},$$
when $P=\mathsf{A}\rho^{\gamma}$ exactly.
\end{Remark}

%\begin{Remark}
By way of precaution, let us verify the conservation of the total mass 
$M$ and the total energy $E$.

First \eqref{1a} can be written as
$$\frac{\partial\rho}{\partial t}=-(\nabla |\rho \vec{v}); $$
Hence
\begin{align*}
\frac{dM}{dt}&=-\int (\nabla |\rho\vec{v})d\mathcal{V}\\
&=-\int_{|\vec{x}|=R}(\rho\vec{v}|\vec{N})d\mathcal{S} =0,
\end{align*}
where the support of $\rho(t,\cdot)$ is supposed to be included
in $\mathfrak{R}=\{|\vec{x}|<R\}$ and $\vec{N}, d\mathcal{S}$ denote
the outer normal vector and the surface area element of
the boundary $\partial\mathfrak{R}=\{|\vec{x}|=R\}$; This shows that
$dM/dt=0$.

Next, \eqref{1a}\eqref{1b} imply
\begin{align*}
\frac{d}{dt}
\int\frac{1}{2}
\rho|\vec{v}|^2d\mathcal{V}
&=\frac{1}{2}
\int\frac{\partial}{\partial t}
\Big(\rho\sum_k(v^k)^2\Big)
d\mathcal{V} \\
&=\int\Big(\rho\sum_kv^k\partial_tv^k+
\frac{1}{2}
\partial_t\rho\sum_k(v^k)^2\Big)\mathcal{V} \\
&=-\int\Big(
\rho\sum_{j,k}v^kv^j\partial_jv^k+
\sum_kv^k\partial_k P+\sum_k\rho v^k\partial_k\Phi+\\
&+\frac{1}{2}\partial_j(\rho v^j)\sum_k(v^k)^2
\Big)d\mathcal{V}\\
&=-\int\Big(
\frac{1}{2}\rho\sum_{j,k}v^j\partial_j(v^k)^2+
\sum_kv^k\partial_k P+\rho\sum_kv^k\partial_k\Phi+\\
&+
\frac{1}{2}\sum_j\partial_j(\rho v^j)\sum_k(v^k)^2\Big)
d\mathcal{V} \\
&=-\int\frac{1}{2}\sum_{j,k}\partial_j
\Big[\rho v^j(v^k)^2\Big] -
\int\sum_k
(v^k\partial_kP+
\rho v^k\partial_k\Phi)
d\mathcal{V}\\
&=-\int(\sum_k\rho v^k\partial_ku+
\sum_k\rho v^k\partial_k\Phi)
d\mathcal{V} \\
&=-\int\Big(
\frac{\partial\Psi}{\partial t}+
\sum_k\rho v^k\partial_k\Phi\Big)
d\mathcal{V} \\
&=-\frac{d}{dt}\int \Psi d\mathcal{V}+
\int\sum_k\partial_k(\rho v^k)\Phi 
d\mathcal{V} \\
&=-\frac{d}{dt}\int\Psi d\mathcal{V}
-\int\frac{\partial\rho}{\partial t}\Phi d\mathcal{V};
\end{align*}
However
\begin{align*}
\int\frac{\partial\rho}{\partial t}\Phi d\mathcal{V}&=-
\mathsf{G}\int\int\partial_t\rho(t,x)\frac{\rho(t,{x}')}{|\vec{x}-\vec{x}'|}d\mathcal{V}({x})
d\mathcal{V}({x}') \\
&=-\frac{\mathsf{G}}{2}\frac{d}{dt}
\int\frac{\rho(t,x)\rho(t,{x}')}{|\vec{x}-\vec{x}'|}d\mathcal{V}({x})d\mathcal{V}({x}') \\
&=\frac{1}{2}\frac{d}{dt}\int\rho\Phi d\mathcal{V};
\end{align*}
This shows $dE/dt=0$. Here we put $\vec{x}=(x^1,x^2,x^3)^{\top}$ and
$\vec{v}=(v^1, v^2, v^3)^{\top}$, and 
denote $\displaystyle \partial_t=\frac{\partial}{\partial t},
\partial_j=\frac{\partial}{\partial x^j} $, while
$k,j$ run $1,2,3$.\\

Moreover the angular momentum 
\begin{equation}
\vec{J}:=\int \vec{x}\times(\rho\vec{v})d\mathcal{V}({x})
\end{equation}
is constant with respect to $t$ along any compactly supported classical solution.
Here
$$\vec{x}\times\vec{v}=
\begin{bmatrix}
x^2v^3-x^3v^2 \\
x^3v^1-x^1v^3 \\
x^1v^2-x^2v^1
\end{bmatrix}
\qquad \mbox{for}\quad \vec{v}=
\begin{bmatrix}
v^1 \\
v^2 \\
v^3
\end{bmatrix}.$$

Let us show it. Note that (\ref{1b}) can be written, under (\ref{1a}), as
\begin{subequations}
\begin{align}
&\frac{\partial}{\partial t}(\rho v^1)+
\frac{\partial}{\partial x^1}(\rho(v^1)^2)
+\frac{\partial}{\partial x^2}(\rho v^1v^2)+
\frac{\partial}{\partial x^3}(\rho v^1v^3)+\frac{\partial P}{\partial x^1}=
-\rho\frac{\partial\Phi}{\partial x^1} \\
&\frac{\partial}{\partial t}(\rho v^2)+
\frac{\partial}{\partial x^1}(\rho v^2v^1)
+\frac{\partial}{\partial x^2}(\rho (v^2)^2)+
\frac{\partial}{\partial x^3}(\rho v^2v^3)+\frac{\partial P}{\partial x^2}=
-\rho\frac{\partial\Phi}{\partial x^2} \\
&\frac{\partial}{\partial t}(\rho v^3)+
\frac{\partial}{\partial x^1}(\rho v^1v^3)
+\frac{\partial}{\partial x^2}(\rho v^2v^3)+
\frac{\partial}{\partial x^3}(\rho (v^3)^2)+\frac{\partial P}{\partial x^3}=
-\rho\frac{\partial\Phi}{\partial x^3} 
\end{align}
\end{subequations}
Then we have
$$\frac{dJ^1}{dt}=-\int(\vec{x}\times\nabla\Phi(t,x))^1\rho (t,x)d\mathcal{V}({x}),
$$
and so on, for $J=(J^1,J^2, J^3)^T$. Here
$$\vec{x}\times\nabla\Phi=
\begin{bmatrix}
(\vec{x}\times\nabla\Phi)^1 \\
\\
(\vec{x}\times\nabla\Phi)^2 \\
\\
(\vec{x}\times\nabla\Phi)^3
\end{bmatrix}
=
\begin{bmatrix}
x^2\frac{\partial\Phi}{\partial x^3}-x^3\frac{\partial\Phi}{\partial x^2} \\
\\
x^3\frac{\partial\Phi}{\partial x^1}-x^1\frac{\partial\Phi}{\partial x^3} \\
\\
x^1\frac{\partial\Phi}{\partial x^2}-x^2\frac{\partial\Phi}{\partial x^1} 
\end{bmatrix}
$$
On the other hand, the differentiation of
the Newton potential (\ref{3}) gives
$$\frac{\partial\Phi}{\partial x^3}=\mathsf{G}\int
\frac{x^3-(x^3)'}{|\vec{x}-\vec{x}'|^3}\rho(t, (x^1)',(x^2)', (x^3)')d\mathcal{V}({x}'),$$
and so on. Hence we see
\begin{align*}
\frac{dJ^1}{dt}&=\mathsf{G}\int\int
\frac{-x^2(x^3-(x^3)')+x^3(x^2-(x^2)')}{|\vec{x}-\vec{x}'|^3}
\rho(t,x)\rho(t,{x}')d\mathcal{V}({x})d\mathcal{V}({x}') \\
&=\mathsf{G}\int\int
\frac{x^2(x^3)'}{|\vec{x}-\vec{x}'|^3}
\rho(t,x)\rho(t,{x}')d\mathcal{V}({x})d\mathcal{V}({x}') \\
&-\mathsf{G}\int\int
\frac{x^3(x^2)'}{|\vec{x}-\vec{x}'|^3}
\rho(t,x)\rho(t,{x}')d\mathcal{V}({x})d\mathcal{V}({x}') \\
&=0.
\end{align*}
By the same manner we can show $dJ^2/dt=dJ^3/dt=0$.$\blacksquare$

\section{Axisymmetric solutions}

\subsection{Co-ordinate system $(\varpi, \phi, z)$}

Let $(\varpi, \phi, z)$ be the cylindrical coordinates defined by
\begin{equation}
x^1=\varpi \cos\phi,\quad x^2=\varpi\sin\phi,\quad x^3=z.
\end{equation}
(This somewhat clumsy notation with $\varpi$ is 
historically standard for this problem.)

Note that, while the polar co-ordinates are
$$x^1=r\sin\theta\cos\phi,\quad
x^2=r\sin\theta\sin\phi, \quad x^3=r\cos\theta,
$$
we are taking $\varpi=r\sin\theta=\sqrt{r^2-z^2}$.\\

We have
\begin{subequations}
\begin{align}
&\frac{\partial}{\partial x^1}=\frac{x^1}{\varpi}\frac{\partial}{\partial\varpi}
-\frac{x^2}{\varpi^2}\frac{\partial}{\partial\phi}, \\
&\frac{\partial}{\partial x^2}=\frac{x^2}{\varpi}\frac{\partial}{\partial \varpi}
+\frac{x^1}{\varpi^2}\frac{\partial}{\partial\phi},\\
&\frac{\partial}{\partial x^3}=\frac{\partial}{\partial z},
\end{align}
\end{subequations}
since
\begin{subequations}
\begin{align}
&\frac{\partial}{\partial\varpi}=\frac{x^1}{\varpi}
\frac{\partial}{\partial x^1}+\frac{x^2}{\varpi}
\frac{\partial}{\partial x^2}, \\
&\frac{\partial}{\partial z}=\frac{\partial}{\partial x^3}, \\
&\frac{\partial}{\partial\phi}=-x^2
\frac{\partial}{\partial x^1}+x^1
\frac{\partial}{\partial x^2}.
\end{align}
\end{subequations}\\

\begin{Definition}
A compactly supported solution $\rho, \vec{v}, \Phi$ will be said to be
{\bf axisymmetric} if $\partial\rho/\partial\phi=0,
\partial\Phi/\partial\phi=0$,
that is,
$\rho=\rho(t,\varpi,z),  \Phi=\Phi(t,\varpi, z)$
and if the velocity field $\vec{v}$ is of the form
\begin{equation}
\vec{v}=
\begin{bmatrix}
\displaystyle \frac{V}{\varpi}x^1-\Omega x^2 \\
\\
\displaystyle \frac{V}{\varpi}x^2+\Omega x^1 \\
\\
W
\end{bmatrix}
\end{equation}
or
\begin{equation}
\vec{v}=V\frac{\partial}{\partial\varpi}+W\frac{\partial}{\partial z}+
\Omega\frac{\partial}{\partial \phi}
\end{equation}
with
$V=V(t,\varpi, z), W=W(t,\varpi, z),
\Omega=\Omega(t,\varpi, z)$.
\end{Definition}

Note that
if $\partial\rho/\partial\phi=0$, then the Newton potential
$\Phi$ given by \eqref{3} necessarily satisfies $\partial\Phi/\partial\phi=0$.

Of course        a spherically symmetric solution, for which
$$\rho=\rho(t,r),\quad \frac{V}{\varpi}=\frac{v(t,r)}{r},\quad \Omega =0,\quad
W=\frac{v(t,r)}{r}z
$$
with $r=\sqrt{\varpi^2+z^2}$, is axisymmetric in this sense.\\

Let us derive the equations which govern axisymmetric solutions.\\

First we nte that the following formula is easily verified:
\begin{equation}
\frac{D}{Dt}=\frac{\partial}{\partial t}+V\frac{\partial}{\partial\varpi}+
\Omega\frac{\partial}{\partial\phi}+W\frac{\partial}{\partial z}.
\end{equation}
In fact, we have
$$x^1\frac{\partial}{\partial x^1}+x^2\frac{\partial}{\partial x^2}=
\varpi\frac{\partial}{\partial\varpi}, \qquad
-x^2\frac{\partial}{\partial x^1}+x^1
\frac{\partial}{\partial x^2}=\frac{\partial}{\partial\phi}.
$$\\

Then the equation \eqref{1a} reads
\begin{equation}
\frac{D\rho}{Dt}+\rho\Big(\frac{\partial V}{\partial \varpi}+
\frac{V}{\varpi}+\frac{\partial W}{\partial z}\Big)=0. \label{X2.4}
\end{equation}\\

Using the calculations
\begin{align*}
&\frac{D x^1}{Dt}=\frac{V}{\varpi}x^1-\Omega x^2,\\
&\frac{D x^2}{Dt}=\frac{V}{\varpi} x^2+\Omega x^1, \\
&\frac{D\varpi}{Dt}=V,
\end{align*}
we see that the equation \eqref{1b} reads
\begin{subequations}
\begin{align}
&\rho\Big[
\frac{x^1}{\varpi}\frac{DV}{Dt}-x^2\frac{D\Omega}{Dt}-
2\frac{x^2}{\varpi}V\Omega -x^1\Omega^2\Big]
+\frac{\partial P}{\partial x^1}=-\rho
\frac{\partial\Phi}{\partial x^1}, \label{X2.5a}\\
&\rho\Big[
\frac{x^2}{\varpi}\frac{DV}{Dt}+x^1\frac{D\Omega}{Dt}+
2\frac{x^1}{\varpi}V\Omega -x^2\Omega^2\Big]
+\frac{\partial P}{\partial x^2}=-\rho
\frac{\partial\Phi}{\partial x^2}, \label{X2.5b}\\
&\rho\frac{DW}{Dt}+\frac{\partial P}{\partial x^3}=-\rho
\frac{\partial\Phi}{\partial x^3}.\label{X2.5c}
\end{align}
\end{subequations}\\

Taking $\displaystyle \frac{1}{\varpi}(x^1\cdot\mbox{\eqref{X2.5a}}+x^2\cdot
\mbox{\eqref{X2.5b}})$ and
$-x^2\cdot\mbox{\eqref{X2.5a}}+x^1\cdot
\mbox{\eqref{X2.5b}})$, we see that \eqref{X2.5a}$\wedge$\eqref{X2.5b} is equivalent to
\begin{subequations}
\begin{align}
&\rho\Big[\frac{DV}{Dt}-\varpi \Omega^2\Big]+\frac{\partial P}{\partial\varpi}=
-\rho\frac{\partial\Phi}{\partial \varpi}, \label{X2.6a}\\
&\rho\frac{D}{Dt}(\varpi^2\Omega)=0.\label{X2.6b}
\end{align}
\end{subequations}\\

%%%%%%%%%%%%%%%%%%%%%%%%%%%%%%%%%%%%%

%%%%%%%%%%%%%%%%%%%%%%%%%%%%%%%%%%%%%%%%%%%%%%%%

On the other hand te Laplace equation (\ref{1c}) reads
\begin{equation}
\frac{1}{\varpi}\frac{\partial}{\partial\varpi}\Big(\varpi\frac{\partial\Phi}{\partial\varpi}\Big)+
\frac{\partial^2\Phi}{\partial z^2}=4\pi \mathsf{G}\rho,\label{14}
\end{equation}
and the Newton potential (\ref{3}) reads
\begin{equation}
\Phi(\varpi,z)=
-4\pi\mathsf{G}\int_{-\infty}^{+\infty}
\int_0^{\infty}
K_I(\varpi,\varpi', z-z')\rho(\varpi',z')\varpi'd\varpi'dz',
\end{equation}
where
\begin{equation}
K_I(\varpi,\varpi',z-z')=
\frac{1}{\pi}\int_0^{\pi/2}
\frac{d\alpha}{\sqrt{(\varpi-\varpi')^2+(z-z')^2+4\varpi\varpi'\sin^2\alpha}}
\end{equation}

Summing up, the full system which governs axisymmetric solutions is

\noindent (\ref{X2.4})(\ref{X2.6a})(\ref{X2.6b})(\ref{X2.5c})(\ref{14}), that is,
\begin{subequations}
\begin{align}
&\frac{D\rho}{Dt}+\rho\Big(\frac{\partial V}{\partial \varpi}+
\frac{V}{\varpi}+\frac{\partial W}{\partial z}\Big)=0,\label{X2.10a}\\
&\rho\Big[\frac{DV}{Dt}-\varpi \Omega^2\Big]+\frac{\partial P}{\partial\varpi}=
-\rho\frac{\partial\Phi}{\partial \varpi}, \label{X2.10b}\\
&\rho\frac{DW}{Dt}+\frac{\partial P}{\partial z}=-\rho
\frac{\partial\Phi}{\partial z}.\label{X2.10c}\\
&\rho\frac{D}{Dt}(\varpi^2\Omega)=0. \label{X2.10d}
\end{align}
\end{subequations}

Here we note the operator $D/Dt$ acting on functions which are axisymmetric
reduces to
$$\frac{D}{Dt}=\frac{\partial}{\partial t}+V\frac{\partial}{\partial\varpi}+
W\frac{\partial}{\partial z}.
$$\\

%%%%%%%%%%%%%%%%%%%%%%%%%%%%%%%%%%%%%%%%%%%%%%%%%%%

Note that \eqref{X2.10d} is a linear first order partial differential equation of $\Omega$,
provided that $\rho\not=0$ and the components of the velocity fields $V, W$
are known. Therefore, given $V, W$, the equation \eqref{X2.10d} can be solved explicitly
as follows.

Let 
\begin{equation}
\Omega^0(\varpi, z)=\Omega|_{t=0}
\end{equation}
be the initial data. For $t\in [0,T), \varpi >0, |z|<\infty$,
we consider the solution
$\tau\mapsto (\varphi(\tau \  ;\     t,\varpi, z), \psi(\tau\  ;\     t,\varpi,z))$
of the ordinary differential equations
$$
\frac{d\varphi}{d\tau}=V(\tau, \varphi, \psi),\qquad
\frac{d\psi}{d\tau}=W(\tau, \varphi,\psi)
$$
satisfying the initial conditions
$$\varphi(t\  ;\    t,\varpi,z)=\varpi,\qquad
\psi(t\  ;\    t,\varpi,z)=z. $$
Then the solution exists on the time interval $[0,t]$,
provided that $V, W$ are bounded, and $\Omega$ is given by
\begin{equation}
\Omega(t,\varpi,z)=\frac{\phi(0\  ;\    t,\varpi,z)^2}{\varpi^2}\Omega^0(
\varphi(0\  ;\    t,\varpi,z), \psi(0\  ;\    t,\varpi,z)).
\end{equation}

Especially let us note that, if $C$ is an arbitrary constant, then
\begin{equation}
\Omega(t, \varpi, z)=\frac{C}{\varpi^2}
\end{equation}
satisfies the equation \eqref{X2.10b}, whatever
$\rho, V, W$ may be, when $\{$ \eqref{X2.10a},\eqref{X2.10b},
\eqref{X2.10c}, \eqref{14} $\}$ 
turns out to be a closed system for only $\rho, V, W \Phi$.
However, if $C\not=0$, then $-\Omega x^2, \Omega x^1$ are unbounded at the axis $\varpi=0$.\\

Let us calculate the angular momentum for the axisymmetric solution.

We see
\begin{align*}
\vec{x}\times\vec{v}&=
\begin{bmatrix}
x^1 \\
\\
x^2 \\
\\
x^3
\end{bmatrix}
\times
\begin{bmatrix}
\displaystyle \frac{V}{\varpi}x^1-\Omega x^2\\
\\
\displaystyle \frac{V}{\varpi}x^2+\Omega x^1 \\
\\
W
\end{bmatrix}
= \\
&=\begin{bmatrix}
\displaystyle-\frac{V}{\varpi}x^3x^2-\Omega x^3x^1 +Wx^2 \\
\\
\displaystyle \frac{V}{\varpi}x^3x^1-\Omega x^3x^2 -Wx^1 \\
\\
\Omega((x^1)^2+(x^2)^2)
\end{bmatrix}
= \\
&=\begin{bmatrix}
-Vz\sin\phi-\Omega z\varpi\cos\phi+W\varpi\sin\phi \\
\\
Vz\cos\phi -\Omega z\varpi\sin\phi-W\varpi\cos\phi \\
\\
\Omega \varpi^2
\end{bmatrix}.
\end{align*}

Integrating this,we get
\begin{equation}
\vec{J}=
\begin{bmatrix}
0 \\
0\\
J
\end{bmatrix},
\end{equation}
where
\begin{equation}
J=2\pi\int_{-\infty}^{+\infty}
\int_0^{\infty}\rho(t,\varpi,z)\Omega(t,\varpi,z)\varpi^3d\varpi dz.
\end{equation}
In fact,
we note that
$$\int_{-\infty}^{+\infty}\int_0^{2\pi}\int_0^{\infty}
\rho V z\sin\phi\cdot\varpi d\varpi d\phi dz=0,
$$
since
$\displaystyle\int_0^{2\pi}\sin\phi d\phi=0$ and so on.

We know that $J$ is constant with respect to $t$ .\\

\begin{Definition}
Let $\rho=\rho(t,\varpi, z), 
\displaystyle \vec{v}=(\frac{V}{\varpi}x^1-\Omega x^2, \frac{V}{\varpi}x^2+\Omega x^1, W)^T$
with $V=V(t,\varpi, z), \Omega=\Omega(t,\varpi, z), W=W(t,\varpi, z)$
be an axisymmetric solution. This solution is said to be {\bf equatorially
symmetric} if 
\begin{align*}
&\rho(t,\varpi, -z)=\rho(t,\varpi, z),\qquad V(t,\varpi,-z)=V(t,\varpi, z), \\
&\Omega(t,\varpi, -z)=\Omega(t,\varpi, z), \qquad
W(t,\varpi,-z)=-W(t,\varpi, z)
\end{align*}
for $\forall z$. 
\end{Definition}
Then the potential $\Phi$ given by \eqref{3} necessarily satisfies
$$\Phi(t,\varpi, -z)=\Phi(t,\varpi, z)$$
for $\forall z$.

\subsection{Co-ordinate system $(r,  \zeta, \phi)$}

Sometimes instead the co-ordinates $(\varpi,z)$ one uses the co-ordinates $(r,\zeta)$ defined by
\begin{equation}
r=\sqrt{\varpi^2+z^2},\qquad \zeta=\frac{z}{r}=
\frac{z}{\sqrt{\varpi^2+z^2}}.
\end{equation}

 Note that, while the polar co-ordinates are
$$x^1=r\sin\theta\cos\phi,\quad
x^2=r\sin\theta\sin\phi,\quad x^3=r\cos\theta, $$
we are taking $\zeta=\cos\theta$ so that
\begin{equation}
x^1=r\sqrt{1-\zeta^2}\cos\phi, \quad x^2=r\sqrt{1-\zeta^2}\sin\phi, \quad x^3=r\zeta,
\end{equation}
while
\begin{subequations}
\begin{align}
&\frac{\partial}{\partial x^1}=\frac{x^1}{r}\frac{\partial}{\partial r}-
\frac{x^1\zeta}{r^2}\frac{\partial}{\partial \zeta}-
\frac{x^2}{r^2(1-\zeta^2)}\frac{\partial}{\partial\phi}, \\
&\frac{\partial}{\partial x^2}=\frac{x^2}{r}\frac{\partial}{\partial r}-
\frac{x^2\zeta}{r^2}\frac{\partial}{\partial \zeta}+
\frac{x^1}{r^2(1-\zeta^2)}\frac{\partial}{\partial\phi}, \\
&\frac{\partial}{\partial x^3}=\zeta\frac{\partial}{\partial r}+
\frac{1-\zeta^2}{r}\frac{\partial}{\partial \zeta},
\end{align}
\end{subequations}
since
\begin{subequations}
\begin{align}
&\frac{\partial}{\partial r}=\frac{x^1}{r}\frac{\partial}{\partial x^1}+
\frac{x^2}{r}\frac{\partial}{\partial x^2}+\zeta\frac{\partial}{\partial x^3}, \\
&\frac{\partial}{\partial\zeta}=-
\frac{x^1\zeta}{1-\zeta^2}\frac{\partial}{\partial x^1}
-\frac{x^2\zeta}{1-\zeta^2}\frac{\partial}{\partial x^2}+
r\frac{\partial}{\partial x^3}, \\
&\frac{\partial}{\partial \phi}=-x^2\frac{\partial}{\partial x^1}+x^1\frac{\partial}{\partial x^2}.
\end{align}
\end{subequations}\\

Let us introduce the velocity component variables $v,w$ by
\begin{equation}
\vec{v}=
\begin{bmatrix}
\displaystyle\Big(\frac{1}{r}v-\frac{\zeta}{1-\zeta^2}w\Big)x^1-\Omega x^2 \\
\\
\displaystyle\Big(\frac{1}{r}v-\frac{\zeta}{1-\zeta^2}w \Big)x^2+\Omega x^1 \\
\\
\zeta v+rw
\end{bmatrix}. \label{vec}
\end{equation}
\\ 

We are going to derive the equations,, assuming that $v, w, \Omega$ are functions of $t, r, \zeta$, and independent of $\phi$.\\

First we note that
the formula
$$\frac{D}{Dt}=\frac{\partial}{\partial t}+v\frac{\partial}{\partial r}+w
\frac{\partial}{\partial \zeta}+\Omega\frac{\partial}{\partial\phi}
$$
holds. Note the following calculations:
\begin{align*}
&\frac{Dx^1}{Dt}=\frac{x^1}{r}v-\frac{x^1\zeta}{1-\zeta^2}w-x^2\Omega,\qquad 
\frac{Dx^2}{Dt}=\frac{x^2}{r}v-\frac{x^2\zeta}{1-\zeta^2}w+x^1\Omega, \\
&\frac{Dr}{Dt}=v,\qquad \frac{D\zeta}{Dt}=w, \\
&\frac{D}{Dt}\Big(\frac{x^1}{r}\Big)=
-\frac{x^1\zeta}{r(1-\zeta^2)}w-\frac{x^2}{r}\Omega,\\
& \frac{D}{Dt}\Big(\frac{x^2}{r}\Big)=-\frac{x^2\zeta}{r(1-\zeta^2)}w+\frac{x^1}{r}\Omega, \\
&\frac{D}{Dt}\Big(\frac{x^1\zeta}{1-\zeta^2}\Big)=
\frac{x^1\zeta}{r(1-\zeta^2)}v+\frac{x^1}{(1-\zeta^2)^2}w-
\frac{x^2\zeta}{1-\zeta^2}\Omega,\\
&\frac{D}{Dt}\Big(\frac{x^2\zeta}{1-\zeta^2}\Big)
=\frac{x^2\zeta}{r(1-\zeta^2)}v+
\frac{x^2}{(1-\zeta^2)^2}w+\frac{x^1\zeta}{1-\zeta^2}\Omega.
\end{align*}

Using these calculations, we see that \eqref{1b} reads
\begin{subequations}
\begin{align}
&\rho\Big[\frac{x^1}{r}\frac{Dv}{Dt}-\frac{x^1\zeta}{1-\zeta^2}
\frac{Dw}{Dt}-x^2\frac{D\Omega}{Dt}  \nonumber \\
&-\frac{2x^1\zeta}{r(1-\zeta^2)}vw
-\frac{x^1}{(1-\zeta^2)^2}w^2
-\frac{2x^2}{r}v\Omega+\frac{2x^2\zeta}{1-\zeta^2}w\Omega
-x^1\Omega^2\Big] +
\frac{\partial P}{\partial x^1}+\rho\frac{\partial\Phi}{\partial x^1}=0, \label{X2.21a} \\
&\rho\Big[\frac{x^2}{r}\frac{Dv}{Dt}
-\frac{x^2\zeta}{1-\zeta^2}\frac{Dw}{Dt}+x^1\frac{D\Omega}{Dt} \nonumber \\
& -\frac{2x^2\zeta}{r(1-\zeta^2)}vw
-\frac{x^2}{(1-\zeta^2)^2}w^2+
\frac{2x^1}{r}v\Omega 
-\frac{2x^1\zeta}{1-\zeta^2}w\Omega -x^2\Omega\Big] +\frac{\partial P}{\partial x^2}+\rho\frac{\partial\Phi}{\partial x^2}=0, \label{X2.21b} \\
&\rho\Big[\zeta
\frac{Dv}{Dt}+r\frac{Dw}{Dt}+2vw\Big]+\frac{\partial P}{\partial x^3}
+\rho\frac{\partial\Phi}{\partial x^3}=0.
\label{X2.21c}
\end{align}
\end{subequations} \\

Taking
$$
\begin{bmatrix}
\frac{x^1}{r} &\frac{x^2}{r}&\zeta \\
\\
-\frac{x^1\zeta}{r^2} & -\frac{x^2\zeta}{r^2} & \frac{1-\zeta^2}{r} \\
\\
-\frac{x^2}{r^2(1-\zeta^2)} & \frac{x^1}{r^2(1-\zeta^2)} & 0
\end{bmatrix}
\begin{bmatrix}
\eqref{X2.21a} \\
\\
\eqref{X2.21b} \\
\\
\eqref{X2.21c}
\end{bmatrix},
$$
we get
the equivalent system:
\begin{subequations}
\begin{align}
&\rho\Big(\frac{Dv}{Dt}-\frac{r}{1-\zeta^2}w^2-r(1-\zeta^2)\Omega^2\Big)+
\frac{\partial P}{\partial r}+\rho\frac{\partial\Phi}{\partial r}=0, \\
&\rho\Big(\frac{D w}{Dt}+
\frac{2}{r}vw+
\frac{\zeta}{1-\zeta^2}w^2+\zeta(1-\zeta^2)\Omega^2\Big)+
\frac{1-\zeta^2}{r^2}\Big(\frac{\partial P}{\partial \zeta}+\rho
\frac{\partial\Phi}{\partial \zeta}\Big)=0,\\
&\rho\Big(\frac{D\Omega}{Dt}
+
\frac{2}{r}v\Omega 
-\frac{2\zeta}{1-\zeta^2}w\Omega \Big)=0. 
\end{align}
\end{subequations}\\

Summing up, the system of equations which governs the evolution of
$\rho, v, w, \Omega$ is:
\begin{subequations}
\begin{align}
&\frac{D\rho}{Dt}+
\rho\Big(\frac{\partial v}{\partial r}+\frac{2}{r}v+
\frac{\partial w}{\partial\zeta}\Big)=0,
\label{6a}\\
&\rho\Big(\frac{Dv}{Dt}-\frac{r}{1-\zeta^2}w^2-r(1-\zeta^2)\Omega^2\Big)+
\frac{\partial P}{\partial r}+\rho\frac{\partial\Phi}{\partial r}=0, \\
&\rho\Big(\frac{Dw}{Dt}+
\frac{2}{r}vw+\frac{\zeta}{1-\zeta^2}w^2+\zeta(1-\zeta^2)\Omega^2\Big)+
\frac{1-\zeta^2}{r^2}\Big(\frac{\partial P}{\partial \zeta}+\rho
\frac{\partial\Phi}{\partial \zeta}\Big)=0,\\
&\rho\frac{D}{Dt}\Big(
r^2(1-\zeta^2)\Omega\Big)=0. \label{202.21d}
\end{align}
\end{subequations}

Here we note that the operator $D/Dt$ reduces to
$$
\frac{D}{Dt}=\frac{\partial}{\partial t}+
v\frac{\partial}{\partial r}+
w\frac{\partial}{\partial \zeta},
$$
when acting on axisymmetric functions $\rho, v, w,\Omega$ 
of $t,r,\zeta$ which do not depend on $\phi$.\\

Of course, a solution $\rho=\rho(t,r,\zeta), v=v(t, r, \zeta), w=w(t, r, \zeta), \Omega=\Omega(t,r,\zeta),
\Phi=\Phi(t,r,\zeta)$ is said to be {\bf equatorially symmetric} if
\begin{align*}
&\rho(t,r,-\zeta)=\rho(t,r,\zeta),\qquad v(t,r,-\zeta)=v(t,r,\zeta) , \qquad
w(t,r,-\zeta)=-w(t,r,\zeta), \\
& \Omega(t,r,-\zeta)=\Omega(t,r, \zeta), \qquad
\Phi(t,r,-\zeta)=\Phi(t,r,\zeta).
\end{align*}
for $\forall \zeta \in [-1,1]$.\\

We note that \eqref{vec} implies
\begin{align}
 |\vec{v}|^2&=(v^1)^2+(v^2)^2+(v^3)^2) = \nonumber\\
&=
\Big[\Big(\frac{1}{r}v-\frac{\zeta}{1-\zeta^2}w\Big)x^1-\Omega x^2\Big]^2
+\Big[\Big(\frac{1}{r}v-\frac{\zeta}{1-\zeta^2}w \Big)x^2+\Omega x^1\Big]^2+
[\zeta v+rw]^2 = \nonumber \\
&=v^2+\frac{r^2}{1-\zeta^2}w^2+r^2(1-\zeta^2)\Omega^2.
\end{align}\\

\begin{Remark}\label{Rem2}

Note that  \eqref{202.21d} can be integrated, 
if $v, w$ and the initial data
$\Omega^0(r, \zeta)=\Omega\Big|_{t=0}$ are given. The solution $\Omega$ is given by
$$\Omega(t,r,\zeta)=\frac{\varphi(0)^2(1-\psi(0)^2)}{r^2(1-\zeta^2)}
\Omega^0(\phi(0), \psi(0)), $$
where the couple of the functions $\tau\mapsto \varphi(\tau)=\varphi(\tau;t,r,\zeta),
\tau\mapsto\psi(\tau)=\psi(\tau;t,r,\zeta)$ is the solution of
\begin{align*}
&\frac{d\varphi}{d\tau}=v(\tau,\varphi,\psi),\quad
\frac{d\psi}{d\tau}=w(\tau, \varphi,\psi), \\
&\varphi(t;t,r,\zeta)=r,\quad
\psi(t;t,r,\zeta)=\zeta.
\end{align*}
\end{Remark}

As for the Poisson equation and the Newton potential, we have
\begin{equation}
\triangle\Phi=\Big[\frac{1}{r^2}
\frac{\partial}{\partial r}r^2\frac{\partial}{\partial r}+\frac{1}{r^2}
\frac{\partial}{\partial \zeta}
(1-\zeta^2)\frac{\partial}{\partial\zeta}\Big]\Phi=4\pi\mathsf{G}\rho
\end{equation}
and
\begin{equation}
\Phi(t,r,\zeta)=-4\pi\mathsf{G}\int_{-1}^1\int_0^{\infty}
K_{II}(r,\zeta, r', \zeta')\rho(t,r',\zeta')r'^2dr'd\zeta',
\end{equation}
where
\begin{equation}
K_{II}(r,\zeta, r', \zeta')=\frac{1}{4\pi}
\int_0^{2\pi}
\frac{d\beta}{\sqrt{r^2+r'^2-2rr'(\sqrt{1-\zeta^2}\sqrt{1-\zeta'^2}\cos\beta+\zeta\zeta')}}.
\end{equation}\\

\begin{Remark}
When $\rho$ is equatorially symmetric, that is,
 $\rho(r,-\zeta)=\rho(r,\zeta)$, then we have the formula
\begin{align}
&\int\frac{\rho(\vec{x}')}{|\vec{x}-\vec{x}'|}d\mathcal{V}({x}') =4\pi
\int_{-1}^1\int_0^{\infty}
K_{II}(r,\zeta, r',\zeta')\rho(r',\zeta')r'^2dr'd\zeta'= \nonumber \\
&=\int_0^1\int_0^{\infty}\sum_{m=0}^{\infty}
f_{2m}(r,r')P_{2m}(\zeta)P_{2m}(\zeta')\rho(r',\zeta')r'^2dr'd\zeta',
\end{align}
where
\begin{equation}
f_n(r,r'):=
\begin{cases}
\displaystyle\frac{1}{r}\Big(\frac{r'}{r}\Big)^n \quad\mbox{for}\quad r'<r \\
\\
\displaystyle\frac{1}{r'}\Big(\frac{r}{r'}\Big)^n\quad\mbox{for}\quad r<r'
\end{cases}
\end{equation}
and $P_n$ are the Legendre polynomials.
\end{Remark}

%%%%%%%%%%%%%%%%%%%%%%%%%%%%%%%%%%%%%%%

%% SECTION 3 %%%%%%%%%%%%%%%%%%%%%%%%%%%%%%%%%%

%%%%%%%%%%%%%%%%%%%%%%%%%%%%%%%%%%%%%%%%%%%%%

\section{Rotating frame of co-ordinates}

\subsection{Rotating frame with a constant angular velocity}

We are going to formulate the Euler-Poisson equations in a frame of reference which is rotating with angular velocity $\vec{\bar{\Omega}}$ around a point $O$ relative to a Newtonian frame. More precisely, let $(t,x)=(t,x^1,x^2,x^3)$ be a Galilean co-ordinates of an inertial system, and
we introduce the rotating co-ordinate system $(t,y^1,y^2,y^3)$ defined by
\begin{subequations}
\begin{align}
&x^1=(\cos\bar{\Omega} t)y^1-(\sin\bar{\Omega} t)y^2, \\
&x^2=(\sin\bar{\Omega} t)y^1+(\cos\bar{\Omega} t)y^2, \\
&x^3=y^3,
\end{align}
\end{subequations}
or,
\begin{subequations}
\begin{align}
&y^1=(\cos\bar{\Omega} t)x^1+(\sin\bar{\Omega} t)x^2, \\
&y^2=-(\sin\bar{\Omega} t)x^1+(\cos\bar{\Omega} t)x^2, \\
&y^3=x^3.
\end{align}
\end{subequations}
Here $\bar{\Omega}$ is a constant. Thus the co-ordinate system $(t,y^1,y^2,y^3)$ is rotating around the $x^3$-axis with the angular velocity $\bar{\Omega}$. We write $\vec{\bar{\Omega}}=(0,0,\bar{\Omega})^{\top}$ in the frame $(t,x)$, that is,
\begin{equation}
\vec{\bar{\Omega}}:=\bar{\Omega}\frac{\partial}{\partial x^3}. \label{LG30}
\end{equation}

Let us consider a vector field
\begin{equation}
\vec{X}(t)=\sum_{j=1}^3X^j(t)\frac{\partial}{\partial x^j}=
\sum_{k=1}^3Y^k(t)\frac{\partial}{\partial y^k}.
\end{equation}
Then we have the formula
\begin{equation}
\frac{d}{dt}\vec{X}(t):=\sum_j\frac{dX^j}{dt}\frac{\partial}{\partial x^j}
=\sum_k\frac{dY^k}{dt}\frac{\partial}{\partial y^k}+\vec{\bar{\Omega}}\times\vec{X}(t). \label{LG32}
\end{equation}

In fact, by the definition, we see
\begin{align*}
&X^1=(\cos\bar{\Omega} t)Y^1-(\sin\bar{\Omega} t)Y^2, \\
&X^2=(\sin\bar{\Omega} t)Y^1+(\cos\bar{\Omega} t)Y^2, \\
&X^3=Y^3.
\end{align*}
Therefore
\begin{align*}
\frac{dX^1}{dt}&=(\cos\bar{\Omega} t)\frac{dY^1}{dt}-(\sin\bar{\Omega} t)
\frac{dY^2}{dt}+\bar{\Omega}(-(\sin\bar{\Omega} t)Y^1-(\cos\bar{\Omega} t)Y^2) \\
&=\sum_k\frac{dY^k}{dt}\frac{\partial x^1}{\partial y^k}-\bar{\Omega} X^2, \\
\frac{dX^2}{dt}&=(\sin\bar{\Omega} t)\frac{dY^1}{dt}+(\cos\bar{\Omega} t)
\frac{dY^2}{dt}+\bar{\Omega}
((\cos\bar{\Omega} t)Y^1-(\sin\bar{\Omega} t)Y^2) \\
&=\sum_k\frac{dY^k}{dt}\frac{\partial x^2}{\partial y^k}+\bar{\Omega} X^1, \\
\frac{dX^3}{dt}&=\frac{dY^3}{dt}=\sum_k
\frac{dY^k}{dt}\frac{\partial x^3}{\partial y^k}+0,
\end{align*}
which yields \eqref{LG32}, since \eqref{LG30} yields
\begin{equation}
\vec{\bar{\Omega}}\times\vec{X}=\bar{\Omega}\Big(-X^2\frac{\partial}{\partial x^1}+X^1\frac{\partial}{\partial x^2}\Big).
\end{equation}

Applying the formula \eqref{LG32} to
\begin{equation}
\vec{X}(t)=\vec{x}:=\sum_jx^j\frac{\partial}{\partial x^j}=\sum_ky^k\frac{\partial}{\partial y^k},
\end{equation}
we have 
\begin{equation}
\vec{v}=\sum_jv^j\frac{\partial}{\partial x^j}
:=\frac{d\vec{x}}{dt}=\sum_j\frac{dx^j}{dt}\frac{\partial}{\partial x^j}
=\vec{u}+\vec{\bar{\Omega}}\times\vec{x}, \label{4.35}
\end{equation}
where 
\begin{equation}
\vec{u}=\sum_ku^k\frac{\partial}{\partial y^k}:=\sum_k
\frac{dy^k}{dt}\frac{\partial}{\partial y^k}.
\end{equation}
Note that, if we write \eqref{4.35} at every components, we have
\begin{align}
v^1&=(\cos\bar{\Omega} t)u^1-(\sin\bar{\Omega} t)u^2-\bar{\Omega} x^2, \nonumber \\
v^2&=(\sin\bar{\Omega} t)u^1+(\cos\bar{\Omega} t)u^2+\bar{\Omega} x^1,\nonumber \\
v^3&=u^3. \label{uv}
\end{align}

Moreover, applying the formula \eqref{LG32} to $\vec{v}=d\vec{x}/dt$ and
$\vec{u}$, we see
\begin{align}
\frac{d\vec{v}}{dt}&=\frac{d\vec{u}}{dt}+\vec{\bar{\Omega}}\times
\frac{d\vec{x}}{dt} \nonumber \\
&=\frac{d\vec{u}}{dt}+\vec{\bar{\Omega}}\times(\vec{u}+\vec{\bar{\Omega}}\times
\vec{x}) \nonumber \\
&=\sum_k\frac{du^k}{dt}\frac{\partial}{\partial y^k}+
\vec{\bar{\Omega}}\times\vec{u}+
\vec{\bar{\Omega}}\times(\vec{u}+\vec{\bar{\Omega}}\times
\vec{x}) \nonumber \\
&=\sum_k\frac{du^k}{dt}\frac{\partial}{\partial y^k}+
2\vec{\bar{\Omega}}\times\vec{u}+
\vec{\bar{\Omega}}\times(\vec{\bar{\Omega}}\times\vec{x}).
\end{align}

This can be written as
\begin{equation}
\sum_j\frac{d^2x^j}{dt^2}\frac{\partial}{\partial x^j}=
\sum_k\frac{d^2y^k}{dt^2}\frac{\partial}{\partial y^k}+
2\vec{\bar{\Omega}}\times\Big(\sum_k\frac{dy^k}{dt}\frac{\partial}{\partial y^k}
\Big)+\vec{\bar{\Omega}}\times\Big(\vec{\bar{\Omega}}\times\Big(
\sum_ky^k\frac{\partial}{\partial y^k}\Big)\Big).
\end{equation}
Note that $$\vec{\bar{\Omega}}=\bar{\Omega}\frac{\partial}{\partial x^3}=
\bar{\Omega}\frac{\partial}{\partial y^3}$$
so that $$\vec{\bar{\Omega}}\times\Big(\sum_ky^k
\frac{\partial}{\partial y^k}\Big)
=\bar{\Omega}\Big(-y^2\frac{\partial}{\partial y^1}+y^1\frac{\partial}{\partial y^2}\Big)$$
and so on.

It is well-known that $-2\vec{\bar{\Omega}}\times \vec{u}$ is called the
`{\bf deflecting} or {\bf Coriolis acceleration}', and $-\vec{\bar{\Omega}}\times
(\vec{\bar{\Omega}}\times\vec{x})$ is called the `{\bf centrifugal
acceleration}'. These, multiplied by the mass, are not real forces but merely apparent forces.\\

The observation so far paraphrases mathematically the discussion developed
in pp. 139-140 of the Batchelor's book \cite{Batchelor}.

\subsection{Euler-Poisson equations in the rotating frame}

  Let us consider the Euler-Poisson equations both in the inertial frame $(t,x^1,x^2,x^3)$ and
in the rotating frame $(t,y^1,y^2,y^3)$ considered in the previous subsection.\\
  
The velocity fields in the frame $(t,x)$ is
\begin{equation}
\vec{v}(t,x)=\sum_jv^j(t,x)\frac{\partial}{\partial x^j},
\end{equation}
and
\begin{equation}
v^j=\frac{dx^j}{dt}
\end{equation}
along the stream line $\displaystyle t\mapsto \vec{x}=\sum_jx^j(t)
\frac{\partial}{\partial x^j}$ of the particle passing $(t,x)$.

The equation of continuity \eqref{1a} is 
\begin{equation}
\frac{D\rho}{Dt}+\rho\sum_j\frac{\partial v^j}{\partial x^j}=0 \label{4.41}
\end{equation}
and
the equation of motion \eqref{1b} is
\begin{equation}
\rho\frac{D\vec{v}}{Dt}+\nabla P=-\rho\nabla\Phi,\label{LG42}
\end{equation}
where 
\begin{align}
&\frac{D\rho}{Dt}=
\frac{\partial\rho}{\partial t}+\sum_{\ell}v^{\ell}\frac{\partial\rho}{\partial x^{\ell}}, \label{4.43}\\
&\frac{D\vec{v}}{Dt}=\sum_j
\Big(\frac{\partial v^j}{\partial t}+
\sum_{\ell}v^{\ell}\frac{\partial v^j}{\partial x^{\ell}}\Big)
\frac{\partial}{\partial x^j}.
\end{align}\\

We are going to write down the equation of motion \eqref{LG42} by dint of the rotating frame $(t,y^1,y^2,y^3)$. 

Consider the stream line $\displaystyle t\mapsto \vec{x}(t)=
\sum_jx^j(t)\frac{\partial}{\partial x^j}$ such that $\displaystyle \vec{v}=
\frac{d\vec{x}}{dt}$,
in which sense we can say
$$
\frac{D\vec{v}}{Dt}=\frac{d\vec{v}}{dt}=\frac{d^2\vec{x}}{dt^2}.$$
The formula \eqref{LG32} tells us 
$$\vec{v}=\vec{u}+\vec{\bar{\Omega}}\times\vec{x},$$
where
\begin{equation}
\vec{u}=\sum_ku^k\frac{\partial}{\partial y^k}=
\sum_k\frac{dy^k}{dt}\frac{\partial}{\partial y^k}
\end{equation}
is the velocity field which appears in the frame $(t,y^1,y^2,y^3)$. 
\begin{Remark} Note that the velocity field is {\bf not} a vector field on $\mathbb{R}^3$ in the sense of the theory of differentiable manifolds, that
is, it is {\bf not} invariant with respect to the change of co-ordinates.
\end{Remark}
The formula \eqref{LG32} tells us
$$ \frac{d\vec{v}}{dt}=\frac{d\vec{u}}{dt}+\vec{\bar{\Omega}}\times\vec{v}.
$$ 
But the formula \eqref{LG32} tells us
\begin{align*}\frac{d\vec{u}}{dt}&=\sum_k\frac{du^k}{dt}
\frac{\partial}{\partial y^k}+
\vec{\bar{\Omega}}\times\vec{u} \\
&=\sum_k\Big(\frac{\partial u^k}{\partial t}+
\sum_{\ell} u^{\ell}\frac{\partial u^k}{\partial y^{\ell}}\Big)
\frac{\partial}{\partial y^k}+\vec{\bar{\Omega}}\times\vec{u}.
\end{align*}

Let us write
\begin{equation}
\frac{Du^k}{Dt}=\frac{\partial u^k}{\partial t}+
\sum_{\ell}u^{\ell}\frac{\partial u^k}{\partial y^{\ell}}.
\end{equation}
Then we can write
\begin{equation}
\frac{d\vec{v}}{dt}=
\sum_k\frac{Du^k}{Dt}\frac{\partial}{\partial y^k}+
2(\vec{\bar{\Omega}}\times\vec{u})+\vec{\bar{\Omega}}\times
(\vec{\bar{\Omega}}\times\vec{x})
\end{equation}
with
$$\vec{u}=\sum_ku^k\frac{\partial}{\partial y^k},\qquad
\vec{x}=\sum_k y^k\frac{\partial}{\partial y^k},
$$
and $u^k$'s are considered as functions of $(t,y^1,y^2,y^3)$. Thus the equation of motion \eqref{LG42} can be written
\begin{align}
\rho\Big[
\sum_k\frac{Du^k}{Dt}\frac{\partial}{\partial y^k}+&2\vec{\bar{\Omega}}\times
\Big(\sum_ku^k\frac{\partial}{\partial y^k}\Big)+
\vec{\bar{\Omega}}\times\Big(\vec{\bar{\Omega}}\times
\Big(\sum_k y^k\frac{\partial}{\partial y^k}\Big)\Big)\Big]+\nonumber \\
&+\sum_k\frac{\partial P}{\partial y^k}\frac{\partial}{\partial y^k}=-
\rho\sum_k\frac{\partial\Phi}{\partial y^k}
\frac{\partial}{\partial y^k},
\end{align}
or, more briefly, 
\begin{equation}
\rho\Big[\frac{D\vec{u}}{Dt}+2\vec{\bar{\Omega}}\times\vec{u}+
\vec{\bar{\Omega}}\times(\vec{\bar{\Omega}}\times\vec{x})\Big]+\nabla P=
-\rho\nabla\Phi.\label{LG49}
\end{equation}
Here note that 
$$\nabla Q=\sum_j\frac{\partial Q}{\partial x^j}\frac{\partial}{\partial x^j}
=\sum_k\frac{\partial Q}{\partial y^k}\frac{\partial}{\partial y^k}$$
for any quantity $Q$. Actually $(y^1,y^2,y^3)$ is an orthonormal
co-ordinate system for each fixed $t$. The equation \eqref{LG49} is nothing but 
\cite[p.295, (13)]{LyndenBO}. (Here we are writing $P, \Phi$ instead of 
$p, -\psi$ of \cite{LyndenBO}, and neglecting other external force
$\mathbf{F}$.)\\

Note that it can be easily verified that the vector product performs in the same manner in the frame $(t,y^1,y^2,y^3)$, that is,
\begin{align}
&\frac{\partial}{\partial y^1}\times\frac{\partial}{\partial y^2}=
-\frac{\partial}{\partial y^2}\times\frac{\partial}{\partial y^1}=
\frac{\partial}{\partial y^3}, \nonumber  \\
&\frac{\partial}{\partial y^2}\times\frac{\partial}{\partial y^3}=
-\frac{\partial}{\partial y^3}\times\frac{\partial}{\partial y^2}=
\frac{\partial}{\partial y^1}, \nonumber  \\
&\frac{\partial}{\partial y^3}\times\frac{\partial}{\partial y^1}=
-\frac{\partial}{\partial y^1}\times\frac{\partial}{\partial y^3}=
\frac{\partial}{\partial y^2}, \nonumber  \\
&\frac{\partial}{\partial y^1}\times\frac{\partial}{\partial y^1}=
\frac{\partial}{\partial y^2}\times\frac{\partial}{\partial y^2}=
\frac{\partial}{\partial y^3}\times\frac{\partial}{\partial y^3}=0.
\end{align}\\

On the other hand, by tedious calculations, it can be verified that
\begin{equation}
\vec{\bar{\Omega}}\times(\vec{\bar{\Omega}}\times\vec{x})=-
\frac{1}{2}\nabla( |\vec{\bar{\Omega}}\times\vec{x}|^2).\label{LG50}
\end{equation}
Here
$$|\vec{\bar{\Omega}}\times\vec{x}|^2=
(\bar{\Omega}^2x^3-\bar{\Omega}^3x^2)^2+
(\bar{\Omega}^3x^1-\bar{\Omega}^1x^3)^2+
(\bar{\Omega}^1x^2-\bar{\Omega}^2x^1)^2$$
for $\vec{\bar{\Omega}}=\sum_j\bar{\Omega}^j\frac{\partial}{\partial x^j},
\vec{x}=\sum_jx^j\frac{\partial}{\partial x^j}$. Using this identity 
\eqref{LG50}, the equation \eqref{LG49} can be written as
\begin{equation}
\rho\Big(\frac{D\vec{u}}{Dt}+2\vec{\bar{\Omega}}\times\vec{u} \Big)+\nabla P=
-\rho\nabla\Big(\Phi-
\frac{1}{2}|\vec{\bar{\Omega}}\times\vec{x}|^2\Big).
\end{equation}
This is nothing but \cite[p.500, (1)]{Lebovitz}. 

Note that this formulation of the equation of motion, in which the centrifugal acceleration is considered as a part of the potential, can be traced back to the paper
\cite{Bjerknes}  by V. Bjerknes in 1929. See \cite[p.11, (A)]{Bjerknes}.\\

On the other hand, by a direct calculation, it can be verified that the equation of continuity 
\eqref{4.41}, which can be written as
$$\frac{D\rho}{Dt}+\rho(\nabla|\vec{v})=0,$$
is written down in the rotating frame as the apparently same form
\begin{equation}
\frac{D\rho}{Dt}+\rho(\nabla|\vec{u})=0. \label{4.52}
\end{equation}
Here we are denoting
\begin{equation}
\frac{D\rho}{Dt}=\frac{\partial\rho}{\partial t}+\sum_ku^k\frac{\partial\rho}{\partial y^k} \label{4.53}
\end{equation}
and
\begin{equation}
(\nabla|\vec{u})=\sum_k\frac{\partial u^k}{\partial y^k} \label{4.54}
\end{equation}

More precisely speaking, \eqref{4.53} should be written as 
$$\frac{D\rho}{Dt}=\Big(\frac{\partial\rho}{\partial t}\Big)_y+\sum_ku^k
\frac{\partial \rho}{\partial y^k}, $$
where the symbol $(\partial/\partial t)_y$ emphasizes that the partial differentiation with respect to $t$ is done keeping $y$ constant. In this symbol, \eqref{4.43} means
$$\frac{D\rho}{Dt}=\Big(\frac{\partial\rho}{\partial t}\Big)_x+\sum_{\ell}v^{\ell}\frac{\partial\rho}{\partial x^{\ell}}.
$$
Note that 
$$\Big(\frac{\partial\rho}{\partial t}\Big)_x\not=
\Big(\frac{\partial\rho}{\partial t}\Big)_y,\qquad
\sum_jv^j\frac{\partial\rho}{\partial x^j}\not=
\sum_ku^k\frac{\partial\rho}{\partial y^k} $$
generally if $\bar{\Omega} \not= 0$.

However, as for \eqref{4.54}, it can be verified by \eqref{uv} that
\begin{equation}
\sum_j\frac{\partial v^j}{\partial x^j}=
\sum_k\frac{\partial u^k}{\partial y^k}, \label{4.55}
\end{equation}
although it is not obvious.
We might write the left-hand side of \eqref{4.55} as $(\nabla_x|\vec{v})$
 and the right-hand side as $(\nabla_y|\vec{u})$,
and the equality as 
$$(\nabla_x|\vec{v})=(\nabla_y|\vec{u}).$$\\

As for the Poisson equation and the Newton potential the expression in the frame $(t,y^1,y^2,y^3)$ is the same as those in $(t,x^1,x^2,x^3)$, since the co-ordinate transformation
$(x^1,x^2,x^3)\mapsto (y^1,y^2,y^3)$ is orthonormal for each fixed $t$. 

In fact the determinant of the Jacobian matrix
$$\begin{bmatrix}
\frac{\partial x^1}{\partial y^1} & \frac{\partial x^1}{\partial y^2}&
\frac{\partial x^1}{\partial y^3} \\
\\
\frac{\partial x^2}{\partial y^1} & \frac{\partial x^2}{\partial y^2} &
\frac{\partial x^2}{\partial y^3} \\
\\
\frac{\partial x^3}{\partial y^1} & \frac{\partial x^3}{\partial y^2} &
\frac{\partial x^3}{\partial y^3}
\end{bmatrix}
=
\begin{bmatrix}
\cos\bar{\Omega} t & -\sin\bar{\Omega} t & 0 \\
\\
\sin\bar{\Omega} t & \cos\bar{\Omega} t & 0 \\
\\
0 & 0 & 1
\end{bmatrix}
$$
is equal to 1 so that
$$d\mathcal{V}=dx^1dx^2dx^3=dy^1dy^2dy^3, $$
and it can be verified easily that
\begin{align*}
&|\vec{x}|^2=(x^1)^2+(x^2)^2+(x^3)^2=(y^1)^2+(y^2)^2+(y^3)^2, \\
&\triangle \Phi=\sum_j\Big(\frac{\partial}{\partial x^j}\Big)^2\Phi
=\sum_k\Big(\frac{\partial}{\partial y^k}\Big)^2\Phi.
\end{align*}
%%%%%%%%%%%%%%%%%%%%%%%%%%%%%%%%%%%%%%%%%%%%%%%%%

\subsection{Axisymmetric solutions in the rotating frame}

Let us observe the Euler-Poisson equations
in the rotating frame considered in the preceding subsections. {\bf However let us write
$(t,x^1,x^2,x^3)$ instead of $(t,y^1,y^2,y^3)$, and $\vec{v}=(v^1,v^2,v^3)$ instead of
$\vec{u}=(u^1,u^2,u^3)$}. Therefore the Euler-Poisson equations in this rotating frame $(t,x^1,x^2,x^3)$ are
\begin{subequations}
\begin{align}
&\frac{D\rho}{Dt}+\rho(\nabla|\vec{v})=0, \label{4.63a}\\
&\rho\Big[\frac{D\vec{v}}{Dt}+2\vec{\bar{\Omega}}\times\vec{v}+
\vec{\bar{\Omega}}\times(\vec{\bar{\Omega}}\times\vec{x})\Big]+\nabla P=-\rho\nabla\Phi, \label{4.63b}\\
&\triangle\Phi=4\pi\mathsf{G}\rho,\label{4.63c}
\end{align}
\end{subequations}
where
\begin{align*}
&\frac{DQ}{Dt}=\frac{\partial Q}{\partial t}+
\sum_{k=1}^3v^k(t,x)\frac{\partial Q}{\partial x^k}\quad\mbox{for}\quad
Q=\rho, v^j, \\
&(\nabla|\vec{v})=\sum_{k=1}^3\frac{\partial v^k}{\partial x^k}, \qquad 
\vec{\bar{\Omega}}=\bar{\Omega}\frac{\partial}{\partial x^3},
\qquad \vec{x}=\sum_{j=1}^3x^j\frac{\partial}{\partial x^j}, \\
&\nabla Q=\sum_{j=1}^3\frac{\partial Q}{\partial x^j}\frac{\partial}{\partial x^j}\quad\mbox{for}\quad Q=P, \Phi,
\end{align*}
and
$$\triangle\Phi=\sum_{j=1}^3\Big(\frac{\partial}{\partial x^j}\Big)^2\Phi.
$$

Supposing that the support of $\rho(t,\cdot)$ is compact, we replace the Poisson equation
by the Newton potential
\begin{equation}
\Phi(t,x)=-\mathsf{G}\int\int\int\frac{\rho(t, x')}{|x-x'|}
dx_1'dx_2'dx_3'.
\end{equation}\\

{\bf (1) }\  We are going to consider axisymmetric solutions by using the co-ordinate
system $(r,\zeta,\phi)$ defined by
\begin{equation}
x^1=r\sqrt{1-\zeta^2}\cos\phi,\quad
x^2=r\sqrt{1-\zeta^2}\sin\phi,\quad
x^3=r\zeta.
\end{equation}
The density $\rho$ and the gravitational potential $\Phi$ is supposed to be functions of
$(t,r,\zeta)$ and the velocity field $\vec{v}=\sum_j v^j\partial/\partial x^j$ is supposed to be of the form
\begin{subequations}
\begin{align}
v^1&=\frac{x^1}{r}v-\frac{x^1\zeta}{1-\zeta^2}w-x^2\omega, \\
v^2&=\frac{x^2}{r}v-\frac{x^2\zeta}{1-\zeta^2}w+x^1\omega, \\
v^3&=\zeta v+rw,
\end{align}
\end{subequations}
where $v,w,\omega$ are functions of $(t,r,\zeta)$.

Then, by the same calculation as that in \S 2.2, we see that the equations 
\eqref{4.63a}\eqref{4.63b} reduce to
\begin{subequations}
\begin{align}
&\frac{D\rho}{Dt}+
\rho\Big(\frac{\partial v}{\partial r}+\frac{2}{r}v+
\frac{\partial w}{\partial\zeta}\Big)=0,
\label{4.106a}\\
&\rho\Big(\frac{Dv}{Dt}-\frac{r}{1-\zeta^2}w^2-r(1-\zeta^2)(\bar{\Omega}+\omega)^2\Big)+
\frac{\partial P}{\partial r}+\rho\frac{\partial\Phi}{\partial r}=0, \label{4.106b}\\
&\rho\Big(\frac{Dw}{Dt}+
\frac{2}{r}vw+\frac{\zeta}{1-\zeta^2}w^2+\zeta(1-\zeta^2)(\bar{\Omega}+\omega)^2\Big)+
\frac{1-\zeta^2}{r^2}\Big(\frac{\partial P}{\partial \zeta}+\rho
\frac{\partial\Phi}{\partial \zeta}\Big)=0,\label{4.106c}\\
&\rho\frac{D}{Dt}\Big(r^2(1-\zeta^2)(\bar{\Omega}+\omega)\Big)=0,\label{4.106d}
\end{align}
\end{subequations}
where
\begin{equation}
\frac{D}{Dt}=\frac{\partial}{\partial t}+v\frac{\partial}{\partial r}+w\frac{\partial}{\partial \zeta}.
\end{equation}
\\

Note that, as observed in Remark \ref{Rem2}, the equation \eqref{4.106d} 
can be integrated immediately once $v, w$ are supposed to be known, provided that the initial data $\overset{\circ}{\omega}(r,\zeta)=\omega(t_0,r,\zeta)$ is given. Then
the equations \eqref{4.106a}\eqref{4.106b}\eqref{4.106c} turn out to be a closed integro-differential evolution system for unknown $\rho, v, w$, provided that $\Phi, 
\overset{\circ}{\omega}$ are given.\\

Of course the Newton potential is given by
\begin{equation}
\Phi(t,r,\zeta)=-4\pi\mathsf{G}\int_{-1}^1\int_0^{\infty}
K_{II}(r,\zeta, r', \zeta')\rho(t,r',\zeta')r'^2dr'd\zeta',\label{5.6}
\end{equation}
where
\begin{equation}
K_{II}(r,\zeta, r', \zeta')=\frac{1}{4\pi}
\int_0^{2\pi}
\frac{d\beta}{\sqrt{r^2+r'^2-2rr'(\sqrt{1-\zeta^2}\sqrt{1-\zeta'^2}\cos\beta+\zeta\zeta')}}.
\end{equation}\\

{\bf (2)} Consider the equations in the co-ordinate system $(\varpi, z, \phi)$ defined by
$$x^1=\varpi\cos\phi,\quad x^2=\varpi\sin\phi, \quad x^3=z.
$$
The velocity field are written as 
\begin{equation}
\vec{v}=\sum v^j\frac{\partial}{\partial x^j}=
V^{\varpi}\frac{\partial}{\partial\varpi}+V^z\frac{\partial}{\partial z}+
V^{\phi}\frac{\partial}{\partial\phi}
\end{equation}
so that
\begin{align}
v^1&=x^1\frac{V^{\varpi}}{\varpi}-x^2V^{\phi} \nonumber \\
v^2&=x^2\frac{V^{\varpi}}{\varpi}+x^1V^{\phi} \nonumber \\
v^3&=V^{z}.
\end{align}
Suppose that $\rho, \Phi, V^{\varpi}, V^z, V^{\phi}$ are functions of $t, \varpi, z$ independent of $\phi$.

Of course
$$\frac{D}{Dt}=\frac{\partial}{\partial t}+
V^{\varpi}\frac{\partial}{\partial\varpi}+V^z\frac{\partial}{\partial z}
+V^{\phi}\frac{\partial}{\partial \phi}.
$$

Since 
\begin{equation}
(\nabla|\vec{v})=\frac{1}{\varpi^2}\frac{\partial}{\partial \varpi}(\varpi^2V^{\varpi})+\frac{\partial V^z}{\partial z},
\end{equation}
\eqref{4.63a} reads
\begin{equation}
\frac{D\rho}{Dt}+\rho
\Big[
\frac{1}{\varpi^2}\frac{\partial}{\partial \varpi}(\varpi^2V^{\varpi})+\frac{\partial V^z}{\partial z}
\Big]=0. \label{0342}
\end{equation}

Since
\begin{align*}
&2\vec{\bar{\Omega}}\times\vec{v}=2\bar{\Omega}\Big(-v^2\frac{\partial}{\partial x^1}+
v^1\frac{\partial}{\partial x^2}\Big), \\
&\vec{\bar{\Omega}}\times(\vec{\bar{\Omega}}\times\vec{x})=\nabla\Big(-\frac{1}{2}
\bar{\Omega}^2\varpi^2\Big),
\end{align*}
\eqref{4.63b} turns out to be the system of equations
\begin{subequations}
\begin{align}
&\rho\Big[\frac{x^1}{\varpi}\frac{DV^{\varpi}}{Dt} -x^2
\frac{DV^{\phi}}{Dt}-2x^2\frac{V^{\varpi}}{\varpi}V^{\phi}-x^1(V^{\phi})^2 + \nonumber \\
&-2\bar{\Omega}
\Big(x^2\frac{V^{\varpi}}{\varpi}+x^1V^{\phi} \Big)+
\frac{\partial}{\partial x^1}\Big(-\frac{1}{2}\bar{\Omega}^2\varpi^2\Big)\Big]+
\frac{\partial P}{\partial x^1}+\rho\frac{\partial\Phi}{\partial x^1}=0, \label{0343a}\\
&\rho\Big[\frac{x^2}{\varpi}\frac{DV^{\varpi}}{Dt}+
x^2\frac{DV^{\phi}}{Dt}+2x^1\frac{V^{\varpi}}{\varpi}V^{\phi}-x^2(V^{\phi})^2+ \nonumber \\
&+2\bar{\Omega}
\Big(x^1\frac{V^{\varpi}}{\varpi}-x^2V^{\phi}
\Big)
+\frac{\partial}{\partial x^2}\Big(-\frac{1}{2}\bar{\Omega}^2\varpi^2\Big)\Big]+
\frac{\partial P}{\partial x^2}+\rho\frac{\partial\Phi}{\partial x^2}=0, \label{0343b}\\
&\rho\frac{DV^z}{Dt}+\frac{\partial P}{\partial x^3}+\rho
\frac{\partial\Phi}{\partial x^3}=0. \label{0343c}
\end{align}
\end{subequations}\\

When $\rho=\rho(\varpi, z), V^{\varpi}=V^z=0, V^{\phi}=V^{\phi}(\varpi, z)$, the
equation \eqref{0342} is satisfied and the set of equations
\eqref{0343a},\eqref{0343b},\eqref{0343c} is reduced to
\begin{subequations}\begin{align}
-\rho x^1(V^{\phi}+\bar{\Omega})^2+\frac{\partial P}{\partial x^1}+
\rho\frac{\partial\Phi}{\partial x^1}&=0, \label{0344a}\\
-\rho x^2(V^{\phi}+\bar{\Omega})^2
+\frac{\partial P}{\partial x^2}+
\rho \frac{\partial\Phi}{\partial x^2}&=0, \label{0344b}\\
\frac{\partial P}{\partial x^3}+\rho \frac{\partial\Phi}{\partial x^3}&=0, \label{0344c}
\end{align}
\end{subequations}
where we read
$$\frac{\partial}{\partial x^j}=\frac{x^j}{\varpi}\frac{\partial}{\partial\varpi} \quad
\mbox{for}\quad  j=1,2,$$
or the set of equations \eqref{0344a},\eqref{0344b},\eqref{0344c} is equivalent to the set
of two equations:
\begin{subequations}\begin{align}
-\rho \varpi(V^{\phi}+\bar{\Omega})^2+\frac{\partial P}{\partial \varpi}+
\rho\frac{\partial\Phi}{\partial \varpi}&=0, \\
\frac{\partial P}{\partial z}+
\rho \frac{\partial\Phi}{\partial z}&=0.
\end{align}
\end{subequations}

%%%%%%%%%%%%%%%%%%%%%%%%%%%%%%%%%%%%%%%%%%%

\section{Stationary solutions}

\subsection{Axisymmetric stationary solutions}

Let us use the co-ordinates $(t,\varpi,z)$.

\begin{Definition}
An axisymmetric solution is said to be {\bf stationary} if $\rho=\rho(\varpi,
z), V=W=0, \Omega=\Omega(\varpi, z), \Phi=\Phi(\varpi, z)$.
\end{Definition}

The equations are reduced to
\begin{subequations}
\begin{align}
-\rho \Omega^2\varpi +\frac{\partial P}{\partial\varpi}
+\rho\frac{\partial\Phi}{\partial\varpi}&=0, \label{P3.1a}\\
\frac{\partial P}{\partial z}+\rho\frac{\partial\Phi}{\partial z}&=0, \label{P3.1b}\\
\triangle\Phi&=4\pi \mathsf{G}\rho, \label{P3.1c}
\end{align}
\end{subequations}
where
\begin{equation}
\triangle\Phi=
\frac{1}{\varpi}\frac{\partial}{\partial\varpi}\Big(
\varpi\frac{\partial\Phi}{\partial\varpi}\Big)+
\frac{\partial^2\Phi}{\partial z^2}.
\end{equation}\\

In our barotropic case, that is, when $P=P(\rho)$ is a given function of $\rho$, we have $\partial\Omega/\partial z=0$, that, is,
$\Omega=\Omega(\varpi)$.

 In fact 
$$\frac{\partial}{\partial z}\Big[\frac{1}{\rho}\times\mbox{(\ref{P3.1a})}\Big] -
\frac{\partial}{\partial \varpi}\Big[\frac{1}{\rho}\times\mbox{(\ref{P3.1b})\Big]}=0$$
implies $\partial(\varpi\Omega^2)/\partial z=0$, keeping in mind that $du=dP/\rho$
so that
 $$\frac{\partial}{\partial z}\Big(\frac{1}{\rho}\frac{\partial P}{\partial\varpi}\Big)=
\frac{\partial}{\partial \varpi}\Big(\frac{1}{\rho}\frac{\partial P}{\partial z}\Big),$$
provided that $u \in C^2$.\\

If $\Omega$ is a nonzero constant, then the stationary solution is called
a {\bf solid rotation}, and if $\Omega$ is not a constant, then the solution is called a {\bf differential rotation}. \\

Using the variable
$\displaystyle u=\int_0^{\rho}\frac{dP}{\rho}$ on the region where $\rho >0$, the system
(\ref{P3.1a})(\ref{P3.1b})(\ref{P3.1c}) is reduced to
\begin{subequations}
\begin{align}
-\varpi\Omega^2+\frac{\partial u}{\partial\varpi}+
\frac{\partial\Phi}{\partial\varpi}&=0, \label{P3.3a}\\
\frac{\partial u}{\partial z}+\frac{\partial\Phi}{\partial z}&=0, \label{P3.3b}\\
\triangle \Phi&=4\pi \mathsf{G}\rho.
\end{align}
\end{subequations}

Then we have
\begin{equation}
-\triangle u=4\pi \mathsf{G}\rho+f(\varpi),
\end{equation}
where
\begin{equation}
f(\varpi):=-\frac{1}{\varpi}\frac{d}{d\varpi}\Big(\varpi^2\Omega^2\Big).
\end{equation}
since $\partial\Omega/\partial z=0$ identically,
that is, since $\Omega=\Omega(\varpi)$.

And (\ref{P3.3a})(\ref{P3.3b}) imply
\begin{equation}
-\int_0^{\varpi}\Omega(\varpi')^2\varpi'd\varpi'+u+\Phi=\mbox{Constant}.
\end{equation}

Note that, if $\Omega$ is a constant $\bar{\Omega}$, then $f=-2\bar{\Omega}^2$, and
we have
\begin{equation}
-\triangle u=4\pi\mathsf{G}\rho-2\bar{\Omega}^2,
\end{equation}
and
\begin{equation}
u+\Phi=\frac{\varpi^2}{2}\bar{\Omega}^2 +\mbox{Const.}
\end{equation}\\

%\subsection{}

On the other hand, if we use the co-ordinates $(t,r,\zeta)$ defined by
\begin{equation}
r=\sqrt{\varpi^2+z^2},\qquad \zeta=\frac{z}{r}=
\frac{z}{\sqrt{\varpi^2+z^2}},
\end{equation}
the equations turn out to be
\begin{subequations}
\begin{align}
-\rho(1-\zeta^2)r\Omega^2+\frac{\partial  P}{\partial r}
+\rho\frac{\partial\Phi}{\partial r}&=0\label{4.1a} \\
\rho\zeta r^2\Omega^2+\frac{\partial P}{\partial\zeta}+
\rho\frac{\partial\Phi}{\partial \zeta}&=0 \label{4.1b} \\
\triangle \Phi &
=4\pi \mathsf{G}\rho,\label{4.1c}
\end{align}
\end{subequations}
where
$$\triangle \Phi=
\frac{1}{r^2}
\frac{\partial}{\partial r}r^2\frac{\partial\Phi}{\partial r}+\frac{1}{r^2}
\frac{\partial}{\partial \zeta}
(1-\zeta^2)\frac{\partial\Phi}{\partial\zeta}.$$\\

If we use the variable $$u=
\int_0^{\rho}\frac{dP}{\rho}, $$
then on the region $\rho>0$ the system (\ref{4.1a})(\ref{4.1b})(\ref{4.1c}) turns out to be
\begin{subequations}
\begin{align}
-(1-\zeta^2)r\Omega^2+\frac{\partial u}{\partial r}+
\frac{\partial\Phi}{\partial r}&=0, \label{P4.4a}\\
\zeta r^2\Omega^2+\frac{\partial u}{\partial\zeta}+
\frac{\partial\Phi}{\partial\zeta}&=0, \label{P4.4b} \\
\triangle\Phi&=4\pi\mathsf{G}\rho. \label{P4.4c} 
\end{align}
\end{subequations}
Suppose that $\Omega$ is a constant $\bar{\Omega}$. Then we have
\begin{equation}
-\triangle u=4\pi \mathsf{G}\rho-2\bar{\Omega}^2,\label{P4.5}
\end{equation}
and
\begin{equation}
u+\Phi=\frac{1}{2}r^2(1-\zeta^2)\bar{\Omega}^2+\mbox{Const.}. \label{P4.6}
\end{equation}
Clearly we see
that (\ref{P4.4a})$\wedge$(\ref{P4.4b})$\wedge$(\ref{P4.4c})$\Leftrightarrow$
(\ref{P4.4c})$\wedge$(\ref{P4.6}).

\begin{Remark}
Let us suppose that $\Omega$ is a constant $\bar{\Omega}\not= 0$. Then the problem cannot be reduced to the single equation \eqref{P4.5}. In fact, if $P=\mathsf{A}\rho^{\gamma}, \frac{6}{5}<\gamma < 2$, and at least if $\Omega$ is sufficiently small, then \eqref{P4.5},
where $$\rho=\Big[\frac{\gamma-1}{\mathsf{A}\gamma}(u\vee 0)
\Big]^{\frac{1}{\gamma-1}},$$ admits a spherically symmetric solution $u=u_s(r)$ such that
$u_s(0)=1 $ and $u_s(r)>0 \Leftrightarrow 0\leq r < R(\Omega^2)$. Here $R(\Omega^2)$ is a finite number depending on $\Omega^2$. However this $u_s(r)$ cannot give a solution of the 
original problem \eqref{P4.4a}\eqref{P4.4b}. In fact the corresponding gravitational
potential $\Phi=\Phi_s$ is spherically symmetric, therefore
$\displaystyle \frac{\partial u_s}{\partial\zeta}=\frac{\partial\Phi_s}{\partial\zeta}=0$. Then \eqref{P4.4b} would imply $\Omega^2=0$, a contradiction.
\end{Remark}

 %%%%%%%%%%%%%%%%%%%%%%%%%%%%%%
 %%%%%%%%%%%%%%%%%%%%%%%%%%%%%
 %%%%%%%%%%%%%%%%%%%%%%%%%%%%%%

\subsection{Normalization}

Suppose $P=\mathsf{A}\rho^{\gamma}$ exactly and $\displaystyle u=\frac{\mathsf{A}\gamma}{\gamma-1}\rho^{\gamma-1}$.\\

Let us normalize the variables. Taking
$r=\mathsf{a} \bar{r}, u=\beta\bar{u},
\Phi=\beta\bar{\Phi}$, 
%and rewriting $r, u, \Phi$ instead of $\bar{r}, \bar{u}, \bar{\Phi}$, 
we have on the region $u>0$
\begin{subequations}
\begin{align}
-\frac{\mathsf{b}}{2}(1-\zeta^2)\bar{r}+
\frac{\partial \bar{u}}{\partial \bar{r}}&=-\frac{\partial\bar{\Phi}}{\partial \bar{r}}, \label{4.9a} \\
\frac{\mathsf{b}}{2}\zeta \bar{r}^2+\frac{\partial \bar{u}}{\partial\zeta}&=
-\frac{\partial\bar{\Phi}}{\partial\zeta}, \label{4.9b} \\
-\bar{\triangle} \bar{u}=\bar{u}^{\nu}-\mathsf{b},\label{4.9c}
\end{align}
\end{subequations}
with
\begin{equation}
\nu=\frac{1}{\gamma-1},\quad
4\pi \mathsf{G}\Big(\frac{\gamma-1}{\mathsf{A}\gamma}\Big)^{\nu}\mathsf{a}^2
\beta^{\frac{2-\gamma}{\gamma-1}}=1,
\quad 2\bar{\Omega}^2\mathsf{a}^2/\beta=\mathsf{b}.
\end{equation}
In order to fix the idea, we may put
$$ \beta=u_{\mathsf{O}}=\frac{\mathsf{A}\gamma}{\gamma-1}\rho_{\mathsf{O}}^{\gamma-1},$$
$\rho_{\mathsf{O}}$ being the central density, when $\bar{u}(\vec{0})=1$.
Then 
$$\mathsf{a}=\sqrt{\frac{\mathsf{A}\gamma}{4\pi\mathsf{G}(\gamma-1)}}\rho_{\mathsf{O}}^{\frac{\gamma-2}{2}},
\qquad
\mathsf{b} =\frac{\bar{\Omega}^2}{2\pi\mathsf{G}\rho_{\mathsf{O}}},\qquad
\rho(r,\zeta)=\rho_{\mathsf{O}}(\bar{u}(r/\mathsf{a},\mathsf{b})\vee 0)^{\nu},
$$
and the Newton potential is
\begin{equation}
\bar{\Phi}(\bar{r},\zeta)=-
\int_{-1}^1\int_0^{\infty}
K_{II}(\bar{r},\bar{r}',\zeta,\zeta')(\bar{u}(\bar{r}',\zeta')\vee 0)^{\nu}\bar{r}'^2d\bar{r}'d\zeta', \label{Pot}
\end{equation}
where $\bar{u}\vee 0=\max\{\bar{u}, 0\}$, since
\begin{equation}
\bar{\triangle}\bar{\Phi}=(\bar{u}\vee 0)^{\nu}.
\end{equation}

Of course we assume that $\bar{u}(0,0)=\bar{u}(0,\zeta) \forall \zeta \in [-1,1]$, since
it is nothing but $\bar{u}(\vec{0})$. On the other hand, it easy to see the potential
$\bar{\Phi}(\bar{r},\zeta)$ defined by (\ref{Pot}) satisfies $\bar{\Phi}(0,0)=\bar{\Phi}(0,\zeta) \forall \zeta\in [-1.1]$.
In fact we have 
$$\frac{\partial}{\partial \zeta}K_{II}(0,r',\zeta,\zeta')=0\quad\mbox{for}\quad \forall (r',\zeta,\zeta').$$\\

Now (\ref{4.9a})$\wedge$(\ref{4.9b}) is equivalent to
\begin{equation}
\bar{\Phi}+\bar{u}=\frac{\mathsf{b}}{4}(1-\zeta^2)\bar{r}^2+\mbox{Const.}\label{4.13}
\end{equation}
which should hold on the region $\bar{u}>0$. \\

\begin{Remark}
  Let us denote $\bar{u}, \bar{r}$ and so on by
$u, r$ and so on. The equation  (\ref{4.9c}) :
$$-\Big(\frac{1}{r^2}\frac{\partial}{\partial r}
r^2\frac{\partial}{\partial r}+
\frac{1}{r^2}\frac{\partial}{\partial \zeta}(1-\zeta^2)\frac{\partial}{\partial\zeta}\Big)u=
u^{\nu} -\mathsf{b},$$
%where $$\mathsf{K}=4\pi \mathsf{G}
%\Big(\frac{\gamma-1}{\mathsf{A}\gamma}\Big)^{\frac{1}{\gamma-1}}
%\qquad \nu=\frac{1}{\gamma-1}, $$
admits a spherically symmetric solution (near to the Lane-Emden function)
$u=U(r)$ such that $U(r)>0 \Leftrightarrow 0\leq r <R(<+\infty)$,
provided that $\gamma >6/5$ and $\bar{\Omega}^2 \ll 1$,
that is, $\mathsf{b} \ll 1$,
which solves
$$-\frac{1}{r^2}\frac{d}{dr}r^2\frac{dU}{dr}=U^{\nu}-\mathsf{b},
\quad U=1+O(r^2)\quad\mbox{as}\quad r\rightarrow 0.$$
But this does not give a solution to our problem if
$\bar{\Omega}\not=0$, that is, $\mathsf{b}\not=0$. In fact,
the Newton potential 
$$\Phi(r,\zeta)=-\int_{-1}^1
\int_0^R
K_{II}(r,r',\zeta,\zeta')U(r')^{\nu}r'^2dr'd\zeta'$$
associated to this spherically symmetric $U$ is spherically symmetric, that
is, $\partial\Phi/\partial\zeta=0$. Then $\partial U/\partial\zeta=\partial
\Phi/\partial\zeta=0$ would imply $\mathsf{b}=0$ in view of (\ref{4.9b}),
a contradiction.
\end{Remark}

\subsection{Spherically symmetric stationary solutions}

Here let us recall the well known results on spherically
symmetric stationary solutions $\rho=\rho(r), V=W={\Omega}=0, \Phi=\Phi(r)$, when the equation of state is the 
exact $\gamma$-law : $P=\mathsf{\mathsf{A}}\rho^{\gamma}$. See \cite{Chandra1939} and
\cite{JosephL}.

The equations are
\begin{align*}
&\frac{dP}{dr}=-\rho\frac{d\Phi}{dr}, \\
&\frac{1}{r^2}\frac{d}{dr}r^2\frac{d\Phi}{dr}=4\pi \mathsf{G}\rho
\end{align*}
and the Newton potential is
\begin{align*}
&\Phi(r)=- 4\pi\mathsf{G}\int_0^{\infty}K_{III}(r,r')\rho(r')r'^2dr', \\
&K_{III}(r,r')=
\min\Big(\frac{1}{r}, 
\frac{1}{r'}\Big)
\end{align*}
Note that 
$$\Phi(r)=- \frac{\mathsf{G}M}{r} \quad\mbox{for}\quad r\geq R
$$
if $\rho(r)=0$ for $r\geq R$ and
$\displaystyle M=4\pi\int_0^R\rho(r)r^2dr$.

The {\bf Lane-Emden equation of index $\nu$} is
$$-\frac{1}{\xi^2}\frac{d}{d\xi}\xi^2\frac{d\theta}{d\xi}=\theta^{\nu}$$
and the {\bf Lane-Emden function of index $\nu$} is the
solution of this equation which satisfies
$$\theta=1,\quad \frac{d\theta}{d\xi}=0\quad
\mbox{at}\quad \xi=0,$$
which will be denoted by $\theta(\xi\  ;\  \nu)$. If
$1\leq \nu <5$, there is a finite $\xi_1(\nu)$ such that
$\theta(\xi\  ;\  \nu)>0$ for $0\leq \xi<\xi_1(\nu)$ and
$\theta(\xi_1(\nu)\  ;\  \nu)=0$.
The numerical table on \cite[Dover Ed., p.96]{Chandra1939} reads
$$ \begin{array}{r|r|r}
\nu & \xi_1(\nu) & \mu_1(\nu) \\ \hline
1.0 & 3.14159 & 3.14159 \\
1.5 & 3.65375 & 2.71406 \\
2.0 & 4.35287 & 2.41105 \\
2.5 & 5.35528 & 2.18720 \\
3.0 & 6.89685 & 2.01824 \\
3.5 & 9.53581 & 1.89056 \\
4.0 & 14.97155 & 1.79723 \\
4.5 & 31.83646 & 1.73780 \\
5.0& \infty & 1.73205
\end{array}$$

Here
$$\mu_1(\nu):=\int_0^{\xi_1(\nu)}\theta(\xi\  ;  \nu)^{\nu}\xi^2d\xi
=-\xi^2\frac{d\theta(\xi\  ;\  \nu)}{d\xi}\Big|_{\xi=\xi_1(\nu)}.
$$

The equilibrium with the central density $\rho|_{r=0}=\rho_{\mathsf{O}}$
is given by the formula
$$\rho=\rho_{\mathsf{O}}\theta(r/\mathsf{a}\  ;\  \nu)^{\nu},
\quad
u=\frac{\mathsf{A}\gamma}{\gamma-1}\rho_{\mathsf{O}}^{\gamma-1}
\theta(r/\mathsf{a}\  ;\  \nu)$$
with
$$\nu=\frac{1}{\gamma-1},\quad
\mathsf{a}=\sqrt{\frac{\mathsf{A}\gamma}{4\pi \mathsf{G}(\gamma-1)}}
\rho_{\mathsf{O}}^{-\frac{2-\gamma}{2}}.$$
The total mass is
$$M=4\pi\Big(\frac{\mathsf{A}\gamma}{4\pi \mathsf{G}(\gamma-1)}\Big)^{3/2}
\rho_{\mathsf{O}}^{\frac{3\gamma-4}{2}}\mu_1(\nu)$$
and the radius is
$$ R=\sqrt{\frac{\mathsf{A}\gamma}{4\pi \mathsf{G}(\gamma-1)}}
\rho_{\mathsf{O}}^{-\frac{2-\gamma}{2}}\xi_1(\nu).$$

Let $1\leq \nu <5 (\Leftrightarrow 6/5<\gamma \leq 2)$.

If 
$\nu \not=3 (\Leftrightarrow \gamma\not=4/3)$, the equilibrium is uniquely determined for each
given total mass $M$. However if $\nu=3 (\Leftrightarrow \gamma
=4/3)$, the total mass does not depend on the central
density, so infinitely many equilibria correspond  to the
same critical mass. This situation is same for the total
energy. Actually we have
\begin{align*}
E&=\frac{4-3\gamma}{\gamma-1}\int P\cdot 4\pi r^2dr \\
&=\frac{4\gamma-3}{\gamma-1}
4\pi \mathsf{A}\Big(\frac{\mathsf{A}\gamma}{4\pi \mathsf{G}(\gamma-1)}\Big)^{3/2}
\rho_{\mathsf{O}}^{\frac{5\gamma-6}{2}}
\int_0^{\xi_1(\nu)}
\theta(\xi\  ;\  \nu)^{\nu+1}\xi^2d\xi
\end{align*}
vanishes for any central density if $\nu=3 (\Leftrightarrow
\gamma=4/3)$. 

In fact the identity
\begin{align*}
E&=\int\frac{P}{\gamma-1}\cdot 4\pi r^2dr+
\frac{1}{2}
\int \rho\Phi\cdot 4\pi r^2dr \\
&=\frac{4-3\gamma}{\gamma-1}\int P\cdot 4\pi r^2dr
\end{align*}
can be shown as follows:\\

 Let $\displaystyle m=\int_0^r\rho\cdot 4\pi r^2dr$.
Then we see
\begin{align*}
\frac{1}{2}\int\rho\Phi\cdot4\pi r^2dr&=\frac{1}{2}\int_0^R\Phi dm 
=\frac{1}{2}\Phi(R)M-\frac{1}{2}
\int_0^Rm\frac{d\Phi}{dr}dr \\
&=-\frac{1}{2}\frac{\mathsf{G}M}{R}-\frac{1}{2}\int_0^R
\frac{\mathsf{G}m^2}{r^2}dr 
\end{align*}
$$\Big(\mbox{since}\quad r^2\frac{d\Phi}{dr}=4\pi \mathsf{G}\int_0^r
\rho  r^2dr=\mathsf{G}m\Big)$$
\begin{align*}
&=-\frac{1}{2}\int_0^R\frac{2\mathsf{G}m}{r}\frac{dm}{dr}dr 
=-\mathsf{G} \int_0^R4\pi \rho m rdr =
\int_0^R\frac{dP}{dr}4\pi r^3dr 
\end{align*}
$$\Big(\mbox{since}\quad \frac{dP}{dr}=-\rho\frac{d\Phi}{dr}=-
\frac{\mathsf{G}\rho m}{r^2}\Big)$$ 
$$
=-3\int_0^R P\cdot 4\pi r^2dr. 
$$
Thus we have the required identity. $\blacksquare$\\

\section{Lagrangian co-ordinate system and Lagrangian variation }

In this section we follow the Lagrangian variation of the equations which has been adopted by astrophysicists.

\subsection{General definition of Lagrangian co-ordinate systems}

Let $-\infty<t_0<T\leq+\infty$ and let us suppose that the velocity field
$$\vec{v}(t,{x})=(v^1(t,{x}), v^2(t,{x}),
v^3(t,{x}))^{\top}$$
is given as a $C^1$-function of $(t,{x})\in
[t_0,T[\times\mathbb{R}^3$.

Then a co-ordinate system $(\mathfrak{t},\mathfrak{x})$ will be called a {\bf Lagrangian co-ordinate system associated with the velocity field $v$} if 1) $\mathfrak{t}=t$, 2)
$(\mathfrak{t},\mathfrak{x})\mapsto (t,x)$ is a $C^1$-diffeomorphism, and 3)
\begin{equation}
\frac{\partial x^j}{\partial\mathfrak{t}}=
v^j(\mathfrak{t}, x(\mathfrak{t}, \mathfrak{x})),\quad j=1,2,3, t_0\leq \forall \mathfrak{t} <T.
\end{equation}

The Lagrangian co-ordinate system is uniquely specified if the initial conditions
\begin{equation}
x^j(t_0,\mathfrak{x})=f^j(\mathfrak{x}),\quad j=1,2,3
\end{equation}
are given. In particular if 
\begin{equation}
x^j(t_0,\mathfrak{x})=\mathfrak{x},\quad j=1,2,3,
\end{equation}
then we shall call $(\mathfrak{t},\mathfrak{x})$ {\bf standard Lagrangian co-ordinate system starting at $t_0$ associated with the velocity field $\vec{v}$}.

\subsection{Lebovitz' Lagrangian change operator and the convective change operator} 

Let $(t, x)=(t, x^1, x^2, x^3)$ be a co-ordinate system of the space-time, which is not necessarily an inertial Galilean system. Let $-\infty
<t_0 < T\leq +\infty$. Suppose that the 
velocity field $\vec{v}(t,x)$ on $[t_0,T[\times\mathbb{R}^3$ and
another velocity field, say, the unperturbed velocity field
$\bar{\vec{v}}(t,x)$ are given. \\

N. R. Lebovitz \cite {Lebovitz} introduced the quantity
$\xi$, which we shall call the `{\bf Lebovitz' displacement}'  defined as follows: Let 
$t \mapsto x=\varphi(t; s, y)=(\varphi^j(t; s, y) | j=1,2,3) $ and
$t \mapsto x=\bar{\varphi}(t; s, y)=(\bar{\varphi}^j(t; s, y) | j=1,2,3) $
be the perturbed and unperturbed flow of the fluid particle passing
the point $(t,x)=(s,y)$, that is, 
\begin{subequations}
\begin{align}
&\frac{\partial}{\partial t}\varphi^j(t;s,y)=v^j(t,\varphi(t;s,y)), \\
&\varphi^j(s;s,y)=y^j,
\end{align}
\end{subequations}
and
\begin{subequations}
\begin{align}
&\frac{\partial}{\partial t}\bar{\varphi}^j(t;s,y)=\bar{v}^j(t,\bar{\varphi}(t;s,y)),
\label{LG3a} \\
&\bar{\varphi}^j(s;s,y)=y^j,\label{LG3b}
\end{align}
\end{subequations}
where $\vec{v}(t,x)=(v^1, v^2, v^3)^{\top}$ and
$\bar{\vec{v}}(t,x)=(\bar{v}^1,\bar{v}^2, \bar{v}^3)^{\top}$ are the perturbed and unperturbed velocity fields; Then N. R. Lebovitz defined the displacement  $\xi(t,x)=(\xi^1(t,x),\xi^2(t,x),\xi^3(t,x))$ by 
\begin{equation}
\xi^j(t,x):=\varphi^j(t; t_0, \bar{\varphi}(t_0;t,x))-x^j,\qquad j=1,2,3.
\end{equation}\\

{\bf We note that, if $\bar{\vec{v}}=0$ (that is, the unperturbed velocity field everywhere vanishes ), then $
\bar{\varphi}(t;s,y)=y$ (that is, the unperturbed flow is static) and the formula
\begin{equation}
x=\mathfrak{x}+\xi(\mathfrak{t},\mathfrak{x}),\quad t=\mathfrak{t}, \label{LDL}
\end{equation}
defines the standard Lagrangian co-ordinate system starting at $t_0$ associated with the velocity field $u$.}

However, if $\bar{\vec{v}}$ does not identically vanish, then $(\mathfrak{t},\mathfrak{x})$ defined by the formula
\eqref{LDL} is not necessarily a Lagrangian co-ordinate system. In fact
the second term of the right-hand side of
the identity
$$\frac{\partial x^j}{\partial\mathfrak{t}}=
v^j(\mathfrak{t},x)+
\sum_k\Big[
\frac{\partial}{\partial y_k}\varphi^j(\mathfrak{t};t_0,y)
\Big]_{y=\bar{\varphi}(t_0;t,\mathfrak{x})}
\frac{\partial}{\partial\mathfrak{t}}\bar{\varphi}^k(t_0;t,\mathfrak{x})
$$
does not necessarily vanish. 

Of course, even if $u$ does not identically vanish, it always holds that
\begin{equation}
\xi^j(t_0,x)=0,\quad j=1,2,3,\forall x.
\end{equation}\\

Moreover N. R. Lebovitz \cite{Lebovitz} introduced the operator $\Delta$,
which we shall call the ` {\bf Lebovitz' Lagrangian change operator}', defined by
\begin{equation}
\Delta Q(t,x)=Q(t, x+\xi(t,x))-\bar{Q}(t,x),
\end{equation}
where $Q$ and $\bar{Q}$ are the values of a physical quantity in the perturbed and unperturbed flows.

By contrast the {\bf `Eulerian difference operator'} $\delta$ is defined by 
\begin{equation}
\delta Q(t,x)=Q(t,x+\xi(t,x))-\bar{Q}(t,x+\xi(t,x)). \label{Edp}
\end{equation}
This definition is slightly different from that of \cite[p.293, (2)]{LyndenBO} by
$Q(t,x)-\bar{Q}(t,x)$. We prefer the definition by  \eqref{Edp} since it gives the simple formula
\begin{equation}
\Delta Q -\Delta\bar{Q}=\delta Q.
\end{equation}

Clearly when the derivatives of $Q$ are Lipschitz continuous we have
\begin{equation}
\Delta Q=\delta Q +(\xi |\nabla)\bar{Q}+O(|\xi|^2),\label{LL1}
\end{equation}
where 
$$(\xi |\nabla)Q(t,x)=\sum_k\xi^k(t,x)\frac{\partial}{\partial x^k}Q(t,x).$$

\begin{Remark}
Note that in some astrophysical literatures, the Lagrangian change operator $\Delta$ is defined as $\delta +(\xi|\nabla)$. We should be careful to avoid confusions.
\end{Remark} 

On the other hand the so-called `{\bf convective rate of change operator}' $D/Dt$ is defined as usually as follows:
\begin{equation}
\frac{D}{Dt}Q(t,x)=\frac{\partial}{\partial t}Q(t,x)+
\sum_jv^j(t,x)\frac{\partial}{\partial x^j}
Q(t,x),
\end{equation}
when the perturbed velocity field $\vec{v}=(v^1,v^2,v^3)^{\top}$ is given.
In parallel with this $D/Dt$, we write
$$
\overline{\frac{D}{Dt}}=\frac{\partial}{\partial t}+\sum_j
\bar{v}^j(t,x)\frac{\partial}{\partial x^j}.$$
Of course
\begin{equation}
\frac{\partial}{\partial t}Q(t,\varphi(t; s,y))=\frac{DQ}{Dt}(t,\varphi(t;s,y))
\end{equation}
along the flow $x=\varphi(t;s,y)$ for any fixed $s,y$.\\

Now, as for the commutation of the operators $\Delta$ and $D/Dt$, D. Lynden-Bell and J. P. Ostriker claimed the following formula.

\begin{Formula}
For any quantity $Q$ the identity
\begin{equation}
\Delta\frac{D}{Dt}Q=\overline{\frac{D}{Dt}}\Delta Q. \label{LG7}
\end{equation}
holds.
\end{Formula}

See \cite[p.294, (9)]{LyndenBO}. But their  proof described in \cite{LyndenBO} seems to fall back on a kind of too much literary rhetoric. So,
here let us describe a simple and stupidly honest proof of this identity.\\

We are going to prove that
\begin{equation}
\clubsuit :=\overline{\frac{D}{Dt}}\Delta Q(t,x)-
\Delta\frac{D}{Dt}Q(t,x)
\end{equation}
vanishes everywhere.

By the definition we have
\begin{align}
\Delta\frac{D}{Dt}Q(t,x)&=\frac{DQ}{Dt}(t,x+\xi)-\overline{\frac{D}{Dt}}
\bar{Q}(t,x) \nonumber \\
&=\partial_tQ(t,x+\xi)+\sum_jv^j(t,x+\xi)\partial_{z^j}Q(t,z)|_{z=x+\xi} \nonumber \\
&-\partial_t\bar{Q}(t,x)-\sum_j\bar{v}^j(t,x)
\partial_{x^j}\bar{Q}(t,x).
\end{align}

Let us calculate $\overline{D/Dt}\Delta Q$. We have
\begin{align*}
&\overline{\frac{D}{Dt}}\Delta Q(t,x)=\frac{\partial}{\partial t}[\Delta Q(t,x)]+
\sum_j\bar{v}^j(t,x)\frac{\partial}{\partial x^j}[\Delta Q(t,x)] \\
&=\partial_tQ(t,x+\xi)+\sum_j\partial_t\xi^j(t,x)
\partial_{z^j}Q(t,z)|_{z=x+\xi}-\partial_t\bar{Q}(t,x) \\
&+\sum_j\bar{v}^j(t,x)
\frac{\partial}{\partial x^j}[Q(t,x+\xi)]-
\sum_j\bar{v}^j(t,x)\partial_{x^j}\bar{Q}(t,x) \\
&=\partial_tQ(t,x+\xi)+\sum_j\partial_t\xi^j(t,x)
\partial_{z^j}Q(t,z)|_{z=x+\xi}-\partial_t\bar{Q}(t,x) \\
&+\sum_j\bar{v}^j(t,x)
(\delta_j^k+\partial_{x^j}\xi^k(t,x))
\partial_{z^j}Q(t, z)|_{z=x+\xi}
-
\sum_j\bar{v}^j(t,x)\partial_{x^j}\bar{Q}(t,x).
\end{align*}
Therefore it holds that
\begin{align}
\clubsuit &=\sum_j\partial_t\xi^j(t,x)\partial_{z^j}Q(t,z)|_{z=x+\xi}+
\sum_j\bar{v}^j(t,x)
\partial_{z^j}Q(t,z)|_{z=x+\xi} \nonumber \\
&+\sum_{j,k}
\bar{v}^j(t,x)\partial_{x^j}\xi^k(t,x)
\partial_{z^k}Q(t,z)|_{z=x+\xi} \nonumber \\
&-\sum_jv^j(t,x+\xi)\partial_{z^j}Q(t,z)|_{z=x+\xi}.\label{LG10}
\end{align}

On the other hand, differentiating
\begin{equation}
\xi^k(t,x)=\varphi^k(t;t_0,\bar{\varphi}(t_0;t,x))-x^k \label{LG11}
\end{equation}
with respect to $t$, we have
\begin{equation}
\partial_t\xi^k(t,x)=u^k(t,x+\xi)+J_0^k(t,x),
\label{LG12}
\end{equation}
where
\begin{equation}
J_0^k(t.x):=\sum_{\ell}\frac{\partial}{\partial t}
[\bar{\varphi}^{\ell}(t_0; t, x)]
\frac{\partial}{\partial y^{\ell}}[\varphi^k(t;t_0,y)]|_{y=\bar{\varphi}(t_0;t,x)}.
\label{LG13}
\end{equation}
By differentiating \eqref{LG11} with respect to $x^j$, we have
\begin{equation}
\partial_{x^j}\xi^k(t,x)=J_j^k(t,x)-\delta_j^k, \label{LG14}
\end{equation}
where
\begin{equation}
J_j^k(t,x)=
\sum_{\ell}\frac{\partial}{\partial x^j}[\bar{\varphi}^{\ell}(t_0;t,x)]
\frac{\partial}{\partial y^{\ell}}
[\varphi^k(t;t_0,y)]|_{y=\bar{\varphi}(t_0;t,x)}.
\label{LG15}
\end{equation}

Substituting \eqref{LG12} \eqref{LG14} into \eqref{LG10}, we see
\begin{equation}
\clubsuit=\sum_k
\Big(J_0^k(t,x)+\sum_j\bar{v}^j(t,x)
J_j^k(t,x)\Big)
\partial_{z^k}Q(t,z)|_{z=x+\xi},
\label{LG16}
\end{equation}
where $J_0^k$ and $J_j^k$ are given by \eqref{LG13} and \eqref{LG15}.
Let us put
\begin{equation}
C^{\ell}(t,x):=
\frac{\partial}{\partial t}[\bar{\varphi}^{\ell}(t_0;t,x)]+
\sum_j\bar{v}^j(t,x)
\frac{\partial}{\partial x^j}[\bar{\varphi}^{\ell}(t_0; t, x)], \label{LG17}
\end{equation}
which gives
\begin{equation}
J_0^k(t,x)+\sum_j\bar{v}^j(t,x)J_j^k(t,x)=
\sum_{\ell}C^{\ell}(t,x)
\frac{\partial}{\partial y^{\ell}}
[\varphi^k(t;t_0,y)]|_{y=\bar{\varphi}(t_0;t,x)}. \label{LG18}
\end{equation}
Thus in order to prove that $\clubsuit$ vanishes, it is sufficient to prove that $C^{\ell}(t,x)=0$ for $\forall \ell =1,2,3$.
So, we are going to evaluate terms of $C^{\ell}(t,x)$.

Let us recall \eqref{LG3a}\eqref{LG3b} and consider 
\begin{equation}
\bar{\psi}_0^j(\tau, s, y):=\frac{\partial}{\partial s}\bar{\varphi}^j(\tau; s, y), \label{LG19}
\end{equation}
which can evaluate
\begin{equation}
\frac{\partial}{\partial t}[\bar{\varphi}^{\ell}(t_0;t,x)]=
\bar{\psi}_0^{\ell}(\tau, s, y)|_{(\tau, s, y)=(t_0,t,x)} . \label{LG20}
\end{equation} 
The functions $\bar{\psi}_0^j$ is determined by the initial value problem
\begin{subequations}
\begin{align}
&\frac{\partial}{\partial\tau}\bar{\psi}_0^j(\tau,s,y)=
\sum_{\ell}\bar{V}_{\ell}^j(\tau,\varphi(\tau;s,y))
\bar{\psi}_0^{\ell}(\tau, s, y) \label{LG22a} \\
&\bar{\psi}_0^j(s,s,y)=-\bar{v}^j(s,y), \label{LG22b}
\end{align}
\end{subequations}
where
\begin{equation}
\bar{V}_{\ell}^j(t,x)=\frac{\partial}{\partial x^{\ell}}\bar{v}^j(t,x).
\label{LG23}
\end{equation}
Actually \eqref{LG22b} comes from the differentiation of \eqref{LG3b}
with respect to to $s$,
since
$$\partial_{\tau}\bar{\varphi}^j(\tau; s, y)|_{\tau=s}=\bar{v}^j(s,\bar{\varphi}(s;s,y))=\bar{v}^j(s,y).$$

Solving \eqref{LG22a}\eqref{LG22b}, we get
\begin{equation}
\bar{\psi}_0(\tau, s, y)=-\exp\Big[
\int_s^{\tau}\bar{V}(\tau',
\bar{\varphi}(\tau';s,y))d\tau'\Big]\bar{\vec{v}}(s,y),
\label{LG24}
\end{equation}
where $\bar{\psi}_0=(\bar{\psi}_0^j | j=1,2,3)$ is a 3-dimensional vector and
$\bar{V}=(\bar{V}_{\ell}^j | \ell, j=1,2,3)$ is a 3 by 3 matrix. 

Next we consider
\begin{equation}
\bar{\psi}_k^j(t,x):=\frac{\partial}{\partial y^k}\bar{\varphi}^j(\tau; s, y),\label{LG25}
\end{equation}
which can evaluate
\begin{equation}
\frac{\partial}{\partial x^j}[\bar{\varphi}^{\ell}(t_0;t,x)]=
\bar{\psi}_j^{\ell}(\tau,s,y)|_{(\tau,s,y)=(t_0,t,x)}.\label{LG26}
\end{equation}
The functions $\bar{\psi}_k^j$ are determined by the initial value problem
\begin{subequations}
\begin{align}
&\frac{\partial}{\partial t}\bar{\psi}_k^j(\tau, s, y)=
\sum_{\ell}\bar{V}_{\ell}^j(\tau,\bar{\varphi}(\tau; s,y))
\bar{\psi}_k^{\ell}(\tau,s,y),\label{LG28a} \\
&\bar{\psi}_k^j(s,s,y)=\delta_k^j.\label{LG28b}
\end{align}
\end{subequations}
Actually \eqref{LG28b} comes from the differentiation fo \eqref{LG3b}
with respect to $y^k$.
Integrating this, we get
\begin{equation}
\bar{\psi}(\tau,s,y)=
\exp\Big[\int_s^{\tau}\bar{V}(\tau',\bar{\varphi}(\tau';s,y))d\tau'\Big],
\label{LG29}
\end{equation}
where $\bar{\psi}=(\bar{\psi}_k^j | k,j=1,2,3)$ is a 3 by 3 matrix. 

Summing up, we see
\begin{align*}
C^{\ell}(t,x)&=-\Big(\exp\Big[\int_t^{t_0}
\bar{V}(\tau',\bar{\varphi}(\tau';t,x))d\tau'\Big]\bar{\vec{v}}(t,x)\Big)^{\ell} + \\
&+\sum_j\bar{v}^j(t,x)
\Big(\exp\Big[\int_t^{t_0}
\bar{V}(\tau',\bar{\varphi}(\tau';t,x))d\tau'\Big]\bar{\vec{v}}(t,x)\Big)_j^{\ell} \\
&=0
\end{align*}
The proof of the fact  $\clubsuit =0$ is complete,
and we can claim that the identity \eqref{LG7} holds for any quantity $Q$. \\

Next, let us observe the commutation of the operators $\Delta$ and the differentiation 
$\partial/\partial x^j$. We claim:

\begin{Formula}
For any quantity $Q$ it holds that
\begin{equation}
d\Delta Q=(\Delta dQ)J+d\bar{Q}(J-I). \label{LG2.1}
\end{equation}
Here, for quantity $A$, $dA$ denotes the co-variant vector
\begin{equation}
dA=\Big(\frac{\partial A}{\partial x^1}\quad  \frac{\partial A}{\partial x^2} \quad  
\frac{\partial A}{\partial x^3}\Big)=\sum_{j=1}^3
\frac{\partial A}{\partial x^j}dx^j
\end{equation}
and we put
\begin{equation}
J=( J_j^k \  |\   j,k=1,2,3),
\end{equation}
with
\begin{align}
J_j^k(t,x)&=\frac{\partial\xi^k}{\partial x^j}+\delta_j^k= \nonumber \\
&=\sum_{\ell}\frac{\partial}{\partial x^j}
[\bar{\varphi}^{\ell}(t_0;t,x)]\frac{\partial}{\partial y^{\ell}}
[\varphi^k(t;t_0,y)]|_{y=\bar{\varphi}(t_0;t,x)}.
\end{align}
(Recall \eqref{LG14} and \eqref{LG15}.) Therefore \eqref{LG2.1} means
\begin{equation}
\frac{\partial}{\partial x^j}(\Delta Q)=
\sum_k\Big(\Big(\Delta\frac{\partial Q}{\partial x^k}\Big)J_j^k+
\frac{\partial \bar{Q}}{\partial x^k}(J_j^k-\delta_j^k)\Big),\quad j=1,2,3,
\end{equation}
or in other words,
\begin{equation}
\Delta\frac{\partial Q}{\partial x^k}=
\sum_j\Big(\frac{\partial}{\partial x^j}(\Delta Q)\cdot(J^{-1})_k^j+
\frac{\partial\bar{Q}}{\partial x^j}((J^{-1})_k^j-\delta_k^j)\Big).
\end{equation}
\end{Formula}

Proof can be done directly. In fact we have
\begin{align*}
&\Delta \frac{\partial}{\partial x^j}Q(t,x)=\partial_{z^j}Q(t,z)|_{z=x+\xi}-\partial_{x^j}\bar{Q}(t,x), \\
&\Big(\frac{\partial}{\partial x^j}\Delta Q\Big)(t,x)=\partial_{z^j}Q(t,z)|_{z=x+\xi}-\partial_{x^j}\bar{Q}(t,x)+ 
\sum_k\frac{\partial\xi^k}{\partial x^j}(t,x)\partial_{z^k}Q(t,z)|_{z=x+\xi}
\end{align*}
and
$$\partial_{z^k}Q(t,z)|_{z=x+\xi}=(\Delta\partial_{x^k}Q)(t,x)+\partial_{x^k}\bar{Q}(t,x) $$
by the definition. \\

Note that, if $\bar{v}^j=0$ for $\forall j$ and $\bar{\varphi}(t;s,y)=y$ identically, then
$$
J_j^k(t,x)=\frac{\partial}{\partial x^j}\varphi(t;t_0, x) 
=\Big(\exp\Big[\int_{t_0}^t
V(\tau',\varphi(\tau'; t_0,x))d\tau'\Big]\Big)_j^k,
$$
where
$$V(t,x)=(V_{\ell}^j(t,x)\  |\  \ell,j=1,2,3)=
\Big(\frac{\partial v^j}{\partial x^{\ell}}(t,x)\  |\  \ell,j=1,2,3\Big).
$$
Therefore, if $\bar{v}^j=0 \forall j$ and
$|\partial v^j/\partial x^{\ell}|\ll 1\  \forall \ell, j$, then $J\doteqdot I$ and 
\eqref{LG2.1} says $d\Delta Q \doteqdot \Delta dQ$. \\

Finally, the following formula can be easily verified:

\begin{Formula}
For any quantities $Q_1, Q_2$, it holds that
\begin{equation}
\Delta(Q_1\cdot Q_2)=(\Delta Q_1)\cdot\bar{Q}_2 +
(\bar{Q}_1+\Delta Q_1)\cdot (\Delta Q_2).
\end{equation}
\end{Formula}

 \begin{Formula}
For any quantity $Q$ and any smooth function $F$, it holds that
\begin{equation}
\Delta F(Q)=F(\bar{Q}+\Delta Q)-F(\bar{Q})
=\Big(\int_0^1
DF(\bar{Q}+\theta\Delta Q)d\theta\Big)\cdot(\Delta Q).
\end{equation}
\end{Formula}

%%%%%%%%%%%%%%%%%%%%%%%%%%%%%%%%%%%%%%%%%%%%%%%%%%%%%%%%%%%%%

%%%%%%%%%%%%%%%%%%%%%%%%%%%%%%%%%%%%%%%%%
\section{Perturbation of the Euler-Poisson equations in the rotating frame}

Let us observe the perturbation of the Euler-Poisson equations
in the rotating frame considered in Section 3.2. {\bf However let us write
$(t,x^1,x^2,x^3)$ instead of $(t,y^1,y^2,y^3)$, and $\vec{v}=(v^1,v^2,v^3)$ instead of
$\vec{u}=(u^1,u^2,u^3)$}. Therefore the Euler-Poisson equations in this rotating frame $(t,x^1,x^2,x^3)$ are
\begin{subequations}
\begin{align}
&\frac{D\rho}{Dt}+\rho(\nabla|\vec{v})=0, \label{54.63a}\\
&\rho\Big[\frac{D\vec{v}}{Dt}+2\vec{\bar{\Omega}}\times\vec{v}+
\vec{\bar{\Omega}}\times(\vec{\bar{\Omega}}\times\vec{x})\Big]+\nabla P=-\rho\nabla\Phi, \label{54.63b}\\
&\triangle\Phi=4\pi\mathsf{G}\rho,\label{54.63c}
\end{align}
\end{subequations}
where
\begin{align*}
&\frac{DQ}{Dt}=\frac{\partial Q}{\partial t}+
\sum_{k=1}^3v^k(t,x)\frac{\partial Q}{\partial x^k}\quad\mbox{for}\quad
Q=\rho, v^j, \\
&(\nabla|\vec{v})=\sum_{k=1}^3\frac{\partial v^k}{\partial x^k}, \qquad 
\vec{\bar{\Omega}}=\bar{\Omega}\frac{\partial}{\partial x^3},
\qquad \vec{x}=\sum_{j=1}^3x^j\frac{\partial}{\partial x^j}, \\
&\nabla Q=\sum_{j=1}^3\frac{\partial Q}{\partial x^j}\frac{\partial}{\partial x^j}\quad\mbox{for}\quad Q=P, \Phi,
\end{align*}
and
$$\triangle\Phi=\sum_{j=1}^3\Big(\frac{\partial}{\partial x^j}\Big)^2\Phi.
$$

Supposing that the support of $\rho(t,\cdot)$ is compact, we replace the Poisson equation
\eqref{4.63c} by the Newton potential
\begin{equation}
\Phi(t,x)=-\mathsf{G}\int\int\int\frac{\rho(t, x')}{|x-x'|}
dx_1'dx_2'dx_3'.
\end{equation}\\

%%%%%%%%%%%%%%%%%%%%%%%%%%%%%%%%%

\subsection{Perturbation from  static solutions}

{\bf Let us consider the unperturbed motion which is static in this rotating frame}, that is, 
we suppose the following:\\

{\bf (S):} {\it There is a solution 
$\rho=\bar{\rho}=\bar{\rho}(x)$ with $v^j=\bar{v}^j=0\  \forall j,\  \Phi=\bar{\Phi}=\bar{\Phi}(x)$. Of course the equation
\begin{equation}
\bar{\rho}\  \vec{\bar{\Omega}}\times(\vec{\bar{\Omega}}\times\vec{x})+\nabla\bar{P}=
-\bar{\rho}\nabla\bar{\Phi} \label{4.65}
\end{equation}
should hold.
Moreover 
\begin{equation}
\mathfrak{R}:=\{ x \in \mathbb{R}^3\  | \  \bar{\rho}(x) >0\}
\end{equation}
is a simply connected domain including the origin $O$. } \\

In this case, since $\bar{\varphi}(t;t_0,x)=x$ identically, we have
\begin{equation}
x^j+\xi^j(t,x)=\varphi^j(t;t_0,x).\label{4.66}
\end{equation}
So, by the definition, we have
\begin{align}
&\Delta\vec{x}=\sum_j\xi^j(t,x)\frac{\partial}{\partial x^j}, \label{4.67}\\
&\Delta v^j(t,x)=v^j(t,x+\xi(t,x))=\frac{\partial}{\partial t}\xi^j(t,x), \label{4.68}
\end{align}
and Formula 1 implies
\begin{equation}
\Big(\Delta\frac{Dv^j}{Dt}\Big)(t,x)=\frac{\partial^2\xi^j}{\partial t^2}(t,x).\label{4.69}
\end{equation}\\

Let us observe the equation of continuity \eqref{54.63a}. Integration of this equation is
\begin{equation}
\rho(t,x+\xi)=\rho(t_0,x)\exp\Big[
-\int_{t_0}^t\mbox{tr}V(\tau, \varphi(\tau; t_0,x))d\tau\Big],\label{4.73}
\end{equation}
where
\begin{equation}
V(t,x)=\Big(\frac{\partial v^j}{\partial x^k}(t,x)\ \Big|\  k,j=1,2,3\Big),
\quad\mbox{and}\quad \mbox{tr}V=\sum_j\frac{\partial v^j}{\partial x^j}
\end{equation}

%However operation of the Lebovitz' Lagrangian change operator $\Delta$ on
%\eqref{4.63a} yields
%\begin{equation}
%\frac{\partial}{\partial t}(\Delta\rho)+(\bar{\rho}+\Delta\rho)
%\sum_j\frac{\partial}{\partial z^j}[v^j(t,z)]|_{z=x+\xi}=0. \label{4.72}
%\end{equation}
%Here we have used Formula 1, 3.
But we can evaluate
\begin{equation}
\frac{\partial}{\partial z^{\ell}}[v^j(t,z)]\Big|_{z=x+\xi}(=V_{\ell}^j(t,x+\xi))=
\sum_k\Big(\frac{\partial J}{\partial t}\Big)_k^j(J^{-1})_{\ell}^k(t,x),
\end{equation}
where the matrix $J=(J_k^j \   |\   k,j=1,2,3)$ is defined by
\begin{equation}
J_k^j(t,x)=\delta_k^j+\frac{\partial \xi^j}{\partial x^k},
\end{equation}
and $J^{-1}(t,x)$ is the inverse matrix of $J(t,x)$. 
In fact, differentiating \eqref{4.66}, we get
\begin{align*}
&\frac{\partial}{\partial t}\xi^j(t,x)=
\frac{\partial}{\partial t}\varphi^j(t;t_0,x)=v^j(t,x+\xi(t,x)), \\
&\frac{\partial^2}{\partial t\partial x^k}\xi^j(t,x)=
\sum_{\ell}\frac{\partial}{\partial z^{\ell}}v^j(t,z)\Big|_{z=x+\xi}\cdot
\Big(\delta_k^{\ell}+\frac{\partial \xi^{\ell}}{\partial x^k}\Big),
\end{align*}
that is,
$$\frac{\partial}{\partial t}J_k^j=\sum_{\ell}V_{\ell}^j(t,x
+\xi(t,x))J_k^{\ell}.$$

Therefore, since
$$\mbox{tr}V(t,\varphi(t;t_0,x))=\mbox{tr}
\Big(\frac{\partial J}{\partial t}J^{-1}\Big)(t,x)=
\frac{\partial}{\partial t}
\log\det J(t,x)$$
and since $J(t_0,x)=I$ (the unit matrix), we see that \eqref{4.73}
reads
\begin{equation}
\rho(t,x+\xi(t,x))=
\bar{\rho}(x)+\Delta\rho(t,x)=
\overset{\circ}{\rho}(x)
\det J(t,x)^{-1}, \label{4.78}
\end{equation}
where $\overset{\circ}{\rho}(x)=\rho(t_0,x)$ is initial data.\\

Here let us look at the linearized approximation of the equation \eqref{4.78}, that is, {\bf hereafter we shall denote $Q\approxeq 0$ if $Q$ is of the magnitude smaller highly than 
 that of $\Delta\rho, v^j, \xi^j$ and their derivatives.}

Since 
\begin{equation}
 \det J(t,x)^{-1}\approxeq 1-\sum_j\frac{\partial\xi^j}{\partial x^j},
\end{equation}
the linearized approximation of the equation \eqref{4.78} is
\begin{equation}
(\Delta\rho)(t,x)\approxeq
(\Delta\rho)(t_0,x)-\bar{\rho}(x)\sum_j\frac{\partial \xi^j}{\partial x^j}(t,x),\label{LL2}
\end{equation}

This means that, once the evolution of  Lebovitz' displacement $\xi(t,x)$ has been solved, then the evolution of the Lagrangian change $\Delta\rho$ of the density is determined, in the sense of linearized approximation, by the divergence of $\xi(t,x)$ provided that the initial perturbation $\Delta\rho(t_0,x)=\rho(t_0,x)-\bar{\rho}(x)$ is given.

\begin{Remark} 
I do not see why D. Lynden-Bell and J. P. Ostriker can claim
\cite[p.296, (19)]{LyndenBO},
that is,
\begin{equation}
\Delta\rho+\bar{\rho}\sum_j\frac{\partial\xi^j}{\partial x^j}=0, \label{LL3}
\end{equation}
even if the initial perturbation $\Delta\rho(t_0, x)$ does not
identically vanish.
Now, \eqref{LL1} says
\begin{equation}
\Delta Q \approxeq \delta Q+(\xi|\nabla)\bar{Q}. \label{LL4}
\end{equation}
Therefore \eqref{LL2} says
$$\Big[\delta \rho+(\xi|\nabla)\rho\Big]_{t_0}^t+
\bar{\rho}(\nabla|\xi)\approxeq 0,
$$
and, if
$\delta\rho+(\xi|\nabla)\rho \approxeq 0$
at $t=t_0$, then
\begin{equation}
\delta\rho+(\xi|\nabla)\bar{\rho}+
\bar{\rho}(\nabla|\xi)\approxeq 0,\label{LL5}
\end{equation}
or
\begin{equation}
\delta\rho+(\nabla|\bar{\rho}\xi)\approxeq 0,\label{LL6}
\end{equation}
which is nothing but \cite[(10)]{Chandra64}. Therefore we should say that \cite{Chandra64} also assumes tacitly that
$\Delta\rho$ vanishes at $t=t_0$. 
\end{Remark}

Now, operating $\Delta$ on the equation of motion \eqref{4.63b} yields
\begin{align}
(\Delta\rho)\vec{\bar{\Omega}}\times(\vec{\bar{\Omega}}\times\vec{x})+&
(\bar{\rho}+\Delta\rho)\Big[\Delta\frac{D\vec{v}}{Dt}+2\vec{\bar{\Omega}}\times\Delta\vec{v}+
\vec{\bar{\Omega}}\times(\vec{\bar{\Omega}}\times\Delta\vec{x})\Big] + \nonumber\\
&+\Delta\nabla P= -(\Delta\rho)(\nabla\bar{\Phi})-(\bar{\rho}+\Delta\rho)\Delta\nabla\Phi
\end{align}
by Formula 3.

Here we have used  \eqref{4.67}, \eqref{4.68}, \eqref{4.69}, and, of course,
$$\displaystyle \frac{\partial^2\vec{\xi}}{\partial t^2},\quad
\frac{\partial\vec{\xi}}{\partial t}, \quad \vec{\xi}$$ mean
$$\displaystyle \sum_j\frac{\partial^2\xi^j}{\partial t^2}(t,x)\frac{\partial}{\partial x^j},\quad
\sum_j\frac{\partial\xi^j}{\partial t}(t,x)\frac{\partial}{\partial x^j},\quad
\sum_j\xi^j(t,x)\frac{\partial}{\partial x^j}$$
respectively. Using \eqref{4.65}, which says
$$\cdot\vec{\bar{\Omega}}\times (\vec{\bar{\Omega}}\times\vec{x})+
\frac{1}{\bar{\rho}}\nabla\bar{P}=-\nabla\bar{\Phi}$$
on $\mathfrak{R}$, we have the equation of motion:

\begin{equation}
\frac{\partial^2\vec{\xi}}{\partial t^2}+2\vec{\bar{\Omega}}\times\frac{\partial\vec{\xi}}{\partial t}+
\vec{\bar{\Omega}}\times(\vec{\bar{\Omega}}\times\vec{\xi})-
\frac{ \Delta\rho}{\bar{\rho}(\bar{\rho}+\Delta\rho)}\nabla\bar{P}+
\frac{1}{\bar{\rho}+\Delta\rho}\Delta\nabla P+\Delta\nabla\Phi=0, \label{0600}
\end{equation}
which should hold on the domain $[0,T]\times\mathfrak{R}$.

Hereafter we shall analyze the equation \eqref{0600} in two ways: First we assume
 {\bf (A)} and next we forget {\bf (A)}.

\subsubsection{Under the assumption {\bf (A)}}

Let us suppose {\bf (A)}. We are going to analyze the term
\begin{equation}
\spadesuit 1:=-\frac{\Delta\rho}{\bar{\rho}(\bar{\rho}+\Delta\rho)}\nabla\bar{P}
+\frac{1}{\bar{\rho}+\Delta\rho}\Delta\nabla P \label{0601}
\end{equation}
or
\begin{equation}
[\spadesuit 1]^j:=-\frac{\Delta\rho}{\bar{\rho}(\bar{\rho}+\Delta\rho)}
\frac{\partial \bar{P}}{\partial x^j}+
\frac{1}{\bar{\rho}+\Delta\rho}\Delta\frac{\partial P}{\partial x^j},\quad j=1,2,3. \label{0602}
\end{equation}
 In order to do it we use the state variable, enthalpy density, $u$ defined by

\begin{equation}
u=\int_0^{\rho}\frac{dP}{\rho},
\end{equation}
under the assumption {\bf (A)}. 
Then we can write
\begin{equation}
[\spadesuit 1]^j=-\frac{\Delta\rho}{\bar{\rho}+\Delta\rho}
\frac{\partial \bar{u}}{\partial x^j}+
\frac{1}{\bar{\rho}+\Delta\rho}\Delta\Big[(\bar{\rho}+\Delta\rho)
\frac{\partial u}{\partial x^j}\Big]. \label{0604}
\end{equation}
Using Formula 3, we see
\begin{align*}
[\spadesuit 1]^j&=-\frac{\Delta\rho}{\bar{\rho}+\Delta\rho}
\frac{\partial \bar{u}}{\partial x^j} + \frac{1}{\bar{\rho}+\Delta\rho}\Big[\Delta\rho\cdot\frac{\partial\bar{u}}{\partial x^j}+
(\bar{\rho}+\Delta\rho)\cdot\Delta \frac{\partial u}{\partial x^j}\Big] =\\
&=\Delta\frac{\partial u}{\partial x^j}.
\end{align*}
In other words, we have
\begin{equation}
\spadesuit 1 =\Delta\nabla u.
\end{equation}

Next we analyze
\begin{equation}
\spadesuit 2:=\Delta\nabla \Phi,
\end{equation}
or
\begin{equation}
[\spadesuit 2]^j =\Delta\frac{\partial \Phi}{\partial x^j},\quad j=1,2,3.
\end{equation}
We look at
\begin{align}
[\spadesuit 2]^j&=\partial_j\Phi(t, x+\xi(t,x))-\partial_j\bar{\Phi}(x) \nonumber \\
&=[\spadesuit 2(1)]^j+[\spadesuit 2(2)]^j, 
\end{align}
where
\begin{subequations}
\begin{align}
[\spadesuit 2(1)]^j&=\partial_j\Phi(t,x+\xi(t,x))-
\partial_j\bar{\Phi}(x+\xi(t,x)), \\
[\spadesuit 2(2)]^j&=\partial_j\bar{\Phi}(x+\xi(t,x))-\partial_j\bar{\Phi}(x).
\end{align}
\end{subequations}
We have
$$
[\spadesuit 2(1)]^j=
-4\pi\mathsf{G}\Big[
\partial_j\mathcal{K}\rho(t,\cdot) (z)-\partial_j\mathcal{K}\bar{\rho} (z)\Big]
$$
with $z=x+\xi(t,x)$. But
$$
\partial_j\mathcal{K}h(z)=-\frac{1}{4\pi}\int\frac{z_j-z_j'}{|\vec{z}-\vec{z}'|^3}
h(z')dz' $$
and
$$
\rho(t,z')-\bar{\rho}(z')=\delta\rho(t, z')
$$
lead us to
\begin{equation}
[\spadesuit 2(1)]^j=-4\pi\mathsf{G}\partial_j\mathcal{K}\delta\rho(t,\cdot)(x+\xi(t,x)).
\end{equation}
On the other hand, the equation for the stationary solutions
$$
\vec{\bar{\Omega}}\times(\vec{\bar{\Omega}}\times\vec{x})+\nabla\bar{u}+\nabla\bar{\Phi}=0
$$
leads us to
\begin{equation}
[\spadesuit 2(2)]^j=
-\vec{\bar{\Omega}}\times (\vec{\bar{\Omega}}\times\vec{\xi}) -
\Delta\frac{\partial\bar{u}}{\partial x^j}.
\end{equation}
Thus we have
\begin{equation}
[\spadesuit 2]^j=
-\vec{\bar{\Omega}}\times (\vec{\bar{\Omega}}\times\vec{\xi}) -
\Delta\frac{\partial\bar{u}}{\partial x^j} -4\pi\mathsf{G}\partial_j\mathcal{K}\delta\rho(t,\cdot)(x+\xi(t,x)).
\end{equation}

Summing up, the equation \eqref{0600} can be written as
\begin{equation}
\frac{\partial^2\xi^j}{\partial t^2}+\Big[2\vec{\bar{\Omega}}\times
\frac{\partial\vec{\xi}}{\partial t}\Big]^j+\Delta\frac{\partial}{\partial x^j}(u-\bar{u})  -4\pi\mathsf{G}\partial_j\mathcal{K}\delta{\rho}(t,\cdot)(x+\xi(t,x))=0
\end{equation}

Moreover, by Formula 2, we see
\begin{equation}
\Delta\frac{\partial}{\partial x^j}(u-\bar{u})=
\sum_k\frac{\partial}{\partial x^k}\Delta(u-\bar{u}) (J^{-1})_j^k.
\end{equation}
Since
$$
\Delta(u-\bar{u})(t,x)=\delta u(t,x),$$
the equation \eqref{0600} reads
\begin{align}
&\frac{\partial^2\xi^j}{\partial t^2}+\Big[2\vec{\bar{\Omega}}\times
\frac{\partial\vec{\xi}}{\partial t}\Big]^j+
\sum_k\frac{\partial}{\partial x^k}\delta u(t,x)\cdot
(J^{-1}(t,x))_j^k + \nonumber \\
& -4\pi\mathsf{G}\partial_j\mathcal{K}\delta{\rho}(t,\cdot)(x+\xi(t,x))=0.
\label{0616}
\end{align}\\

Now let us derive the linearized approximation of the equation \eqref{0616}, assuming that $\overset{\circ}{\rho}=\bar{\rho}$, when
$$ \delta\rho  \approxeq -(\nabla|\bar{\rho}\vec{\xi}) $$
and
$$
\delta u(t,x)\approxeq \overline{\frac{du}{d\rho}}(\nabla|\bar{\rho}\vec{\xi})).
$$
Since
$$(J^{-1})_j^k=\delta_j^k+O(D_{\vec{x}}\vec{\xi}),$$
we have the linearized approximation of \eqref{0616} as
\begin{subequations}
\begin{align}
\frac{\partial^2\xi^1}{\partial t^2}&-2\bar{\Omega}\frac{\partial \xi^2}{\partial t}&+
\frac{\partial}{\partial x^1}\Big(-\overline{\frac{1}{\rho}\frac{dP}{d\rho}}g +4\pi\mathsf{G}\mathcal{K}g\Big)
\approxeq 0, \label{EKVAa}\\
\frac{\partial^2\xi^2}{\partial t^2}&+2\bar{\Omega}\frac{\partial \xi^1}{\partial t}&+
\frac{\partial}{\partial x^2}\Big(-\overline{\frac{1}{\rho}\frac{dP}{d\rho}}g +4\pi\mathsf{G}\mathcal{K}g\Big)
\approxeq 0, \label{EKVAb}\\
\frac{\partial^2\xi^3}{\partial t^2}&&+
\frac{\partial}{\partial x^3}\Big(-\overline{\frac{1}{\rho}\frac{dP}{d\rho}}g +4\pi\mathsf{G}\mathcal{K}g\Big)
\approxeq 0, \label{EKVAc}
\end{align}
\end{subequations}
where
\begin{equation}
g:= \sum_k\frac{\partial}{\partial x^k}(\bar{\rho}\xi^k).\label{EKVAg}
\end{equation}

\subsubsection{Without {\bf (A)} supposed on perturbations}

Let us forget {\bf (A)} at the moment, that is, we do not use the barotropic assumption that $P$ is a prescribed function of $\rho$. However we preserve the assumption {\bf (A)} on the unperturbed state, that is, $\bar{P}$ is a prescribed function of
$\bar{\rho}$ independent of the place $x$. Moreover, as in the preceding subsection, we
suppose that  $\overset{\circ}{q}=0$, that is, $\Delta\rho $ vanishes at $t=t_0$. \\

We should consider
$$
[\spadesuit 1]^j=-\frac{\Delta\rho}{\bar{\rho}(\bar{\rho}+\Delta\rho)}
\frac{\partial\bar{P}}{\partial x^j}+
\frac{1}{\bar{\rho}+\Delta\rho}
\Delta\frac{\partial P}{\partial x^j}
$$
without using the variable $u$ of the preceding subsection provided {\bf (A)}.

On the other hand we have
$$
[\spadesuit 2(1)]^j=-4\pi\mathsf{G}\partial_j\mathcal{K}\delta\rho(t,\cdot)(x+\xi(t,x))
$$
as before, but  
$$
[\spadesuit 2(2)]^j=\partial_j\bar{\Phi}(x+\xi(t,x))
-\partial_j\bar{\Phi}(x)
$$
turns out to be
\begin{align*}
[\spadesuit 2(2)]^j&=
-[\vec{\bar{\Omega}}\times(\vec{\bar{\Omega}}\times\vec{\xi})]^j
-\Delta\Big(\frac{1}{\bar{\rho}}\frac{\partial\bar{P}}{\partial x^j}\Big) = \\
&=-[\vec{\bar{\Omega}}\times(\vec{\bar{\Omega}}\times\vec{\xi})]^j+
\frac{\Delta\bar{\rho}}{\bar{\rho}(\bar{\rho}+\Delta\rho)}
\frac{\partial\bar{P}}{\partial x^j}
-\frac{1}{\bar{\rho}+\Delta\bar{\rho}}
\Delta\frac{\partial\bar{P}}{\partial x^j}.
\end{align*}
Therefore
\begin{align*}
[\spadesuit 1+\spadesuit 2]^j&=
-[\vec{\bar{\Omega}}\times(\vec{\bar{\Omega}}\times\vec{\xi})]^j + \\
&+
\frac{\Delta\bar{\rho}-\Delta\rho}{(\bar{\rho}+\Delta\rho)(\bar{\rho}+\Delta\bar{\rho})}\frac{\partial\bar{P}}{\partial x^j}
-\frac{1}{\bar{\rho}+\Delta\bar{\rho}}\Delta\frac{\partial\bar{P}}{\partial x^j} + 
\frac{1}{\bar{\rho}+\Delta\rho}\Delta\frac{\partial P}{\partial x^j} +\\
&-4\pi\mathsf{G}\partial_j\mathcal{K}\delta\rho(t,\cdot)(\vec{x}+\xi(t,x)).
\end{align*}
But, using Formula 2, we have
\begin{align*}
\Delta\frac{\partial\bar{P}}{\partial x^j}&=
\sum_k\frac{\partial\Delta\bar{P}}{\partial x^k}(J^{-1})_j^k
-\frac{\partial\bar{P}}{\partial x_j}+
\sum_k\frac{\partial\bar{P}}{\partial x^k}(J^{-1})_j^k, \\
\Delta\frac{\partial{P}}{\partial x^j}&=
\sum_k\frac{\partial\Delta P}{\partial x^k}(J^{-1})_j^k
-\frac{\partial\bar{P}}{\partial x_j}+
\sum_k\frac{\partial\bar{P}}{\partial x^k}(J^{-1})_j^k.
\end{align*}
Noting that
\begin{align*}
(\Delta P-\Delta\bar{P})(t,x)&=\delta P(t,x), \\
(\bar{P}+\Delta\bar{P})(t,x)&=\bar{P}(x+\xi(t,x)),
\end{align*}
we have
\begin{align}
[\spadesuit 1+\spadesuit 2]^j&=
-[\vec{\bar{\Omega}}\times(\vec{\bar{\Omega}}\times\vec{\xi})]^j + \nonumber \\
&+
\frac{\Delta\bar{\rho}-\Delta\rho}{(\bar{\rho}+\Delta\rho)(\bar{\rho}+\Delta\bar{\rho})}
\sum_k\frac{\partial}{\partial x_k}\bar{P}(t,x+\xi(t,x))
\cdot(J^{-1}(t,x))_j^k + \nonumber \\
&+\frac{1}{\bar{\rho}+\Delta\rho}
\sum_k\frac{\partial}{\partial x_k}\delta P(t,x)
\cdot
(J^{-1}(t,x))_j^k +\nonumber \\
&-4\pi\mathsf{G}
\partial_j\mathcal{K}\delta\rho(t,\cdot)(x+\xi(t,x)).
\end{align}
Summing up, we see that the equation \eqref{0600} is reduced to
\begin{align}
&\frac{\partial^2\xi^j}{\partial t^2}+
\Big[
2\vec{\bar{\Omega}}\times
\frac{ \partial\vec{\xi} }{\partial t}\Big]^j+ 
\nonumber \\
&+\frac{\Delta\bar{\rho}-\Delta\rho}{(\bar{\rho}+\Delta\rho)(\bar{\rho}+\Delta\bar{\rho})}
\sum_k\frac{\partial}{\partial x_k}\bar{P}(t,\vec{x}+\xi(t,x))\cdot(J^{-1}(t,x))_j^k + \nonumber \\
&+\frac{1}{\bar{\rho}+\Delta\rho}
\sum_k\frac{\partial}{\partial x_k}\delta P(t,x)\cdot
(J^{-1}(t,x))_j^k +\nonumber \\
&-4\pi\mathsf{G}
\partial_j\mathcal{K}\delta\rho(t,\cdot)(x+\xi(t,x))=0.
\end{align}

Now it is clear that the linearized approximation is
\begin{equation}
\frac{\partial^2\vec{\xi}}{\partial t^2}+
2\vec{\bar{\Omega}}\times\frac{\partial\vec{\xi}}{\partial t}
-\frac{\delta\rho}{\bar{\rho}^2}\nabla\bar{P}
+\frac{1}{\bar{\rho}}\nabla\delta P+
\nabla4\pi\mathsf{G}\mathcal{K}(\nabla|\bar{\rho}\vec{\xi})\approxeq 0.
\end{equation}

\subsection{Perturbation from non-static solutions}

Let us describe a generalization to the case in which the unperturbed solution is not
static, say, $\vec{\bar{v}}\not=\vec{0}$. We suppose\\

{\bf (S-):} {\it There is a solution $\rho=\bar{\rho}(x), \vec{v}=\vec{\bar{v}}(x), \Phi=\bar{\Phi}(x)$ such that 
\begin{equation}
(\nabla | \bar{\rho}\vec{\bar{v}})=0\quad\mbox{and}\quad
\frac{D}{Dt}\vec{\bar{v}}=0,
\end{equation}
which satisfies
\begin{equation}
\bar{\rho}\Big[2\vec{\bar{\Omega}}\times\vec{\bar{v}}+
\vec{\bar{\Omega}}\times(\vec{\bar{\Omega}}\times\vec{x})\Big]+
\nabla \bar{P}+\bar{\rho}\nabla\bar{\Phi}=0.
\end{equation}
}\\

The perturbation of the equation of continuity \eqref{54.63a} gives 
\begin{equation}
\bar{\rho}(x)+\Delta \rho(t,x)=\overset{\circ}{\rho}(x)\det J(t,x)^{-1}
\end{equation}
which is the same to \eqref{4.78}. But the perturbation of the equation of motion 
\eqref{54.63b} is different. Firstly we note that \eqref{4.68} should be replaced by
\begin{equation}
\Big(\Delta\frac{Dv^j}{Dt}\Big)(t,x)=
\frac{\partial^2\xi^j}{\partial t^2}(t,x)+
\sum_k\bar{v}^k(x)\frac{\partial^2\xi^j}{\partial x^k\partial t}(t,x).
\end{equation}
Here we have used Formula 1 and $$\sum_k\bar{v}^k\frac{\partial \bar{v}^j}{\partial x^k}=0.$$
Briefly writing in the vector form, this is 
\begin{equation}
\Delta \frac{D\vec{v}}{Dt}=\frac{\partial^2\vec{\xi}}{\partial t^2}+(\vec{\bar{v}}|\nabla)\frac{\partial\vec{\xi}}{\partial t}.
\end{equation}
Therefore the operation of $\Delta$on \eqref{54.63b} using Formula 3 gives
\begin{align}
&(\bar{\rho}+\Delta\rho)\Big[\frac{\partial^2\vec{\xi}}{\partial t^2}+
(\vec{\bar{v}}|\nabla)\frac{\partial\vec{\xi}}{\partial t}+
2\vec{\bar{\Omega}}\times\frac{\partial\vec{\xi}}{\partial t}+
\vec{\bar{\Omega}}\times(\vec{\bar{\Omega}}\times\vec{\xi})\Big]+ \nonumber \\
&-\bar{\rho}(2\vec{\bar{\Omega}}\times\vec{\bar{v}})+(\Delta\rho)(\vec{\bar{\Omega}}\times
(\vec{\bar{\Omega}}\times\vec{x}))+ \nonumber \\
&+\Delta\nabla P+(\Delta\rho)\nabla\bar{\Phi}+
(\bar{\rho}+\Delta\rho)\Delta\nabla\Phi=0.
\end{align}
Using the equation
$$
\bar{\rho}\Big[2\vec{\bar{\Omega}}\times\vec{\bar{v}}+
\vec{\bar{\Omega}}\times(\vec{\bar{\Omega}}\times\vec{x})\Big]+
\nabla \bar{P}+\bar{\rho}\nabla\bar{\Phi}=0,
$$
we have
\begin{align}
&\frac{\partial^2\vec{\xi}}{\partial t^2}+(\vec{\bar{v}}|\nabla)\frac{\partial\vec{\xi}}{\partial t}+
2\vec{\bar{\Omega}}\times\frac{\partial\vec{\xi}}{\partial t}
+\vec{\bar{\Omega}}\times(\vec{\bar{\Omega}}\times\vec{\xi})
-(2\vec{\bar{\Omega}}\times\vec{\bar{v}}) + \nonumber \\
&-\frac{\Delta\rho}{\bar{\rho}(\bar{\rho}+\Delta\rho)}\nabla\bar{P}+
\frac{1}{\bar{\rho}+\Delta\rho}\Delta\nabla P+
\Delta \nabla\Phi=0. \label{0655}
\end{align}
Naturally this reduces to \eqref{0600}
when $\vec{\bar{v}}=\vec{0}$.

Let us omit the following calculations, which involve nothing new, e.g., under the assumption {\bf (A)}, the equation \eqref{0655} reduces to
\begin{align}
&\frac{\partial^2{\xi}^j}{\partial t^2}+(\vec{\bar{v}}|\nabla)\frac{\partial{\xi}^j}{\partial t}+
\Big[2\vec{\bar{\Omega}}\times\frac{\partial\vec{\xi}}{\partial t}
+\vec{\bar{\Omega}}\times(\vec{\bar{\Omega}}\times\vec{\xi})
-(2\vec{\bar{\Omega}}\times\vec{\bar{v}})\Big]^j + \nonumber \\
&+
\sum_k\frac{\partial}{\partial x^k}\delta u(t,x)\cdot
(J^{-1}(t,x))_j^k  -4\pi\mathsf{G}\partial_j\mathcal{K}\delta{\rho}(t,\cdot)(x+\xi(t,x))=0,
\end{align}
and so on.

\subsection{Non-zero initial perturbations}

Let us go back to the case $\bar{\vec{v}}=\vec{0}$.\\

The perturbation $\xi(t,x)=\Delta x $ should vanish at $t=0$. But there sometime turns out to be necessary to consider the equations  for $\tilde{\xi}=\xi+\overset{\circ}{\xi}$
with given $\overset{\circ}{\xi}=\overset{\circ}{\xi}(x)\not=0$. In order to consider this case, taking $\overset{\circ}{\rho}=\bar{\rho}$ is not convenient. Alternatively we choose 
$\overset{\circ}{\rho}\not=\bar{\rho}$ in a suitable way. In order to treat it, we suppose\\

{\bf (A+):} $$ P=\mathsf{A}\rho^{\gamma},\quad 1<\gamma <2$$\\

\noindent for the sake of simplicity of the argument. Then
$$u=\frac{\mathsf{A}\gamma}{\gamma-1}\rho^{\gamma-1}\quad\mbox{for}\quad \rho >0.$$
But we permit that $u(t,x)$ takes negative values somewhere, and we put
\begin{equation}
\rho=f^{\rho}(u):=\Big(\frac{\gamma-1}{\mathsf{A}\gamma}\Big)^{\frac{1}{\gamma-1}}
(u \vee 0)^{\frac{1}{\gamma-1}}.
\end{equation}
Here and hereafter  we denote $ Q\vee Q'=\max\{ Q, Q'\}$. \\

We suppose that the equilibrium $\bar{\rho}$ satisfies\\

{\it The equilibrium
$$\bar{u}(x)=\frac{\mathsf{A}\gamma}{\gamma-1}\bar{\rho}(x)^{\gamma-1} \quad (x \in \mathfrak{R}:=\{ x |\bar{\rho}(x) >0\} )$$
admits an extension onto $\mathbb{R}^3$ in the class $C^1(\mathbb{R}^3)$ so that
$\bar{\rho}=f^{\rho}(\bar{u})$, that is, $\bar{u}(x)\leq 0 \Leftrightarrow x \not\in \mathfrak{R}$. }\\

Moreover we interpret the Euler equations \eqref{54.63a}, \eqref{54.63b} as
\begin{subequations}
\begin{align}
&\frac{Du}{Dt}+(\gamma-1)u(\nabla|\vec{v})=0, \label{6.3.1a} \\
&\frac{D\vec{v}}{Dt}+2\vec{\bar{\Omega}}\times\vec{v}+
\vec{\bar{\Omega}}\times(\vec{\bar{\Omega}}\times\vec{x})+
\nabla u+ \nabla\Phi=0 \label{6.3.1b}
\end{align}
\end{subequations}
which should hold on the whole space $[0,T[\times\mathbb{R}^3$. Of course the set 
$u=\bar{u}(x), \vec{v}=0$ enjoys \eqref{6.3.1a}\eqref{6.3.1b}. 

Given the initial data $\overset{\circ}{u}(x)=u(0,x)$ on $\mathbb{R}^3$, we have
\begin{equation}
u(t, x+\xi(t,x))=\bar{u}(x)+\Delta u(t,x)=\overset{\circ}{u}(x)\det J(t,x)^{-(\gamma-1)},
\end{equation}
where
\begin{equation}
J(t,x)=(J_k^j(t,x))_{j,k}, \qquad J_k^i=\delta_k^j+\frac{\partial \xi^j}{\partial x^k}.
\end{equation}\\

Given $\overset{\circ}{\xi}(x)$, we consider an initial data $\overset{\circ}{u}$ such that 
\begin{equation}
\overset{\circ}{u}=(1-(\overset{\circ}{\xi}|\nabla))\bar{u}-
(\gamma-1)\bar{u}(\nabla|\overset{\circ}{\xi})
=
\bar{u}-\overline{\frac{du}{d\rho}}(\nabla|\bar{\rho}\overset{\circ}{\xi})
\quad\mbox{on}\quad \mathfrak{R}.
\end{equation}

Namely, we consider the evolution of $(u, \vec{v})$ governed by the equations 
\eqref{6.3.1a}, \eqref{6.3.1b} with 
\begin{equation}
\Phi=-4\pi\mathsf{G}\mathcal{K}\rho,
\quad
\rho=f^{\rho}(u):=\Big(\frac{\gamma-1}{\mathsf{A}\gamma}\Big)^{\frac{1}{\gamma-1}}
(u \vee 0)^{\frac{1}{\gamma-1}}
\end{equation}
under the initial conditions $u|_{t=0}=\overset{\circ}{u}, \vec{v}|_{t=0}=\overset{\circ}{\vec{v}}$.
The stream function $\varphi(t, x)$ is defined by this velocity field $\vec{v}$ as
the solutions of the initial value problem of the set of ordinary differential equations
$$\partial_t\varphi(t,x)=\vec{v}(t,\varphi(t,x)),\quad \varphi(0,x)=x, $$
and 
$\xi$ is defined by
$ x+\xi(t,x)=\varphi(t, x)$. We should remark that, in this sense, $\xi$ depends upon $\overset{\circ}{u}, \overset{\circ}{\vec{v}}$, and that it can be different with that for $\overset{\circ}{\rho}=\bar{\rho}$, say, $\overset{\circ}{u}=\bar{u}$, which was supposed in the preceding discussions.

Then we have
\begin{align*}
\delta u &=u-\bar{u} \\
&=u-\overset{\circ}{u}+\overset{\circ}{u}-\bar{u} \\
&=\overset{\circ}{u}(\det J^{-(\gamma-1)}-1)-
(\overset{\circ}{\xi}|\nabla)\bar{u}-
(\gamma-1)\bar{u}(\nabla|\overset{\circ}{\xi}).
\end{align*}
So, the linearized approximation of
$$\sum_k\frac{\partial}{\partial x^k}\delta u\cdot (J^{-1})_j^k=
\sum_k\frac{\partial}{\partial x^k}(u-\overset{\circ}{u})\cdot (J^{-1})_j^k+
\sum_k\frac{\partial}{\partial x^k}(\overset{\circ}{u}-\bar{u})\cdot (J^{-1})_j^k
$$
turns out to be
\begin{align*}
&\frac{\partial}{\partial x^j}\Big[
-(\gamma-1)\overset{\circ}{u}(\nabla|\xi)-(\xi|\nabla)\overset{\circ}{u}
-(\gamma-1)\bar{u}(\nabla|\overset{\circ}{\xi})-(\overset{\circ}{\xi}|\nabla)\bar{u} \Big]\\
&\approxeq
\frac{\partial}{\partial x^j}\Big[
-(\gamma-1)\bar{u}(\nabla|\xi)-(\xi|\nabla)\bar{u}
-(\gamma-1)\bar{u}(\nabla|\overset{\circ}{\xi})-(\overset{\circ}{\xi}|\nabla)\bar{u} \Big]\\
&=\frac{\partial}{\partial x^j}\Big[
-(\gamma-1)\bar{u}(\nabla|\xi+\overset{\circ}{\xi})-(\xi+\overset{\circ}{\xi}|\nabla)\bar{u} \Big]\\
&=\frac{\partial}{\partial x^j}\tilde{G}, 
\end{align*}
where 
$$
\tilde{G}:=-(\gamma-1)\bar{u}(\nabla|\xi+\overset{\circ}{\xi})-(\xi+\overset{\circ}{\xi}|\nabla)\bar{u} =-\overline{\frac{du}{d\rho}}(\nabla|\bar{\rho}(\xi+\overset{\circ}{\xi})).
$$

In this sense the system of the linearized equations for $\tilde{\xi}=\xi+\overset{\circ}{\xi}$ is the same to \eqref{EKVAa}\eqref{EKVAb}\eqref{EKVAc}\eqref{EKVAg} with $\xi$ replaced by $\tilde{\xi}$.\\

However we should note that the assumed initial distribution of $u$:
$$\overset{\circ}{u}=
\bar{u}\Big[1-(\gamma-1)(\nabla|\overset{\circ}{\xi})-\Big(\overset{\circ}{\xi}\Big|\frac{\nabla\bar{u}}{\bar{u}}\Big)\Big] $$
can take negative values in $\mathfrak{R}$ near the boundary $\partial\mathfrak{R}$, even if
$\overset{\circ}{\xi}, (\nabla|\overset{\circ}{\xi})$ are very small, since
$\displaystyle \frac{\nabla\bar{u}}{\bar{u}}=
(\gamma-1)\frac{\nabla\bar{\rho}}{\bar{\rho}}$ usually diverges at the vacuum boundary $\partial\mathfrak{R}$. Thus, when $\overset{\circ}{\xi} \not=0$ near the boundary, it will happens that we must consider the  initial density $\displaystyle \overset{\circ}{\rho}
:=f^{\rho}(\overset{\circ}{u})=
\Big(\frac{\gamma-1}{\mathsf{A}\gamma}\Big)^{\frac{1}{\gamma-1}}
(\overset{\circ}{u} \vee 0)^{\frac{1}{\gamma-1}}$ which vanishes somewhere in $\mathfrak{R}$, say, vacuum at the initial time, near the boundary, while we imagine $u$ with negative values enjoys the evolution on the whole space according to the extended Euler equations
\eqref{6.3.1a}\eqref{6.3.1b}. Here we also note that, since
$$\overset{\circ}{u}=
\bar{u}-\overline{\frac{du}{d\rho}}(\nabla|\bar{\rho}\overset{\circ}{\xi})
=\bar{u}
\Big(1-(\gamma-1)\frac{1}{\bar{\rho}}(\nabla|\bar{\rho}\overset{\circ}{\xi})\Big),
$$
we have
$$\overset{\circ}{\rho}=\bar{\rho}\Big[
\Big(1-(\gamma-1)\frac{1}{\bar{\rho}}(\nabla|\bar{\rho}\overset{\circ}{\xi})\Big)
\vee 0\Big]^{\frac{1}{\gamma-1}}.
$$

\vspace{15mm}

%%%%%%%%%%%%%%%%%%%%%%%%%%%%%%%%%%%%%%%%%%%%%%

{\bf\Large Acknowledgment}\\

   I would like to express sincere thanks 
to Professor Juhi Jang who carefully read the first draft of this note
and pointed out mistakes and insufficient expositions in it.
The first draft of this note was prepared for the discussion with her during my stay
at Korea Institute for Advanced Study during July 11 -15, 2016.
I would like to express sincere thanks to KIAS and Professor Juhi Jang
for the full financial support and hospitality for this stay. 
The retouch and revision of this note was done during my stay at the Department of Mathematics of National University of Singapore during
February 8-21, 2017 thanks to the kind invitation by Professor 
Shih-Hsien Yu. I would like to express sincere thanks to the Department of Mathematics of NUS and Professor Yu for the full financial support and hospitality for this stay.

\end{document}